\RequirePackage{fix-cm}
\documentclass{svjour3}                     
\smartqed  
\usepackage[dvipdfmx]{graphicx,color}
\usepackage{bm}
\usepackage[numbers]{natbib}
\bibpunct[:]{(}{)}{,}{a}{}{,}
\usepackage{textcomp} 
\usepackage{amsmath,amssymb,amsfonts} 
\usepackage{setspace}
\usepackage{algorithmic}
\usepackage{algorithm}



\usepackage[switch,pagewise]{lineno}
\newcommand*\patchAmsMathEnvironmentForLineno[1]{
    \expandafter\let\csname old#1\expandafter\endcsname\csname #1\endcsname
    \expandafter\let\csname oldend#1\expandafter\endcsname\csname end#1\endcsname
    \renewenvironment{#1}
    {\linenomath\csname old#1\endcsname}
    {\csname oldend#1\endcsname\endlinenomath}}
\newcommand*\patchBothAmsMathEnvironmentsForLineno[1]{
    \patchAmsMathEnvironmentForLineno{#1}
    \patchAmsMathEnvironmentForLineno{#1*}}
\AtBeginDocument{
\patchBothAmsMathEnvironmentsForLineno{equation}
\patchBothAmsMathEnvironmentsForLineno{align}
\patchBothAmsMathEnvironmentsForLineno{flalign}
\patchBothAmsMathEnvironmentsForLineno{alignat}
\patchBothAmsMathEnvironmentsForLineno{gather}
\patchBothAmsMathEnvironmentsForLineno{multline}
}

\usepackage{url}
\usepackage[toc,title,page]{appendix}

\usepackage[top=30truemm,bottom=30truemm,left=25truemm,right=25truemm]{geometry}
\begin{document}

\title{Self-exciting negative binomial distribution process and critical properties of intensity distribution
}


\author{Kotaro Sakuraba         \and
        Wataru Kurebayashi      \and
        Masato Hisakado         \and
        Shintaro Mori             
}


\institute{
K.Sakuraba, S.Mori \at
Department of Mathematics and Physics,
Graduate School of Science and Technology, 
Hirosaki University \\
Bunkyo-cho 3, Hirosaki, Aomori 036-8561, Japan
\\
W. Kurebayashi \at
Department of Mechanical Science and Engineering,
Graduate School of Science and Technology, 
Hirosaki University \\
Bunkyo-cho 3, Hirosaki, Aomori 036-8561, Japan
\\
M. Hisakado \at
Nomura Holdings, Inc., Otemachi 2-2-2, Chiyoda-ku, Tokyo 100-8130, Japan
}

\date{Received: date / Accepted: date}

\maketitle




\begin{abstract}
\label{abst}
We study the continuous time limit of a self-exciting negative binomial process and discuss the critical properties of its intensity distribution.
In this limit, the process transforms into a marked Hawkes process. 
The probability mass function of the marks has a parameter $\omega$, and the process reduces to a "pure" Hawkes process in the limit $\omega\to 0$. We investigate the Lagrange--Charpit equations for the master equations of the marked Hawkes process
in the Laplace representation
close to its critical point and extend the previous findings on the
power-law scaling of the probability density function (PDF) of intensities in the intermediate asymptotic regime  to the case where the memory kernel is the superposition of an arbitrary 
finite number of exponentials.  
We develop an efficient sampling method for the marked Hawkes process
based on the time-rescaling theorem and verify the 
power-law exponents.
\end{abstract}

\section{Introduction}
\label{intro}
Recently, high-frequency financial data have become available, and 
many studies have been devoted to calibrating models of market micro-structure\citep{Bacry:2015}.
Initially, many studies adopted discrete time models 
that captured trade dynamics at regular time intervals or through a sequence of discrete time events, such as trades.
Subsequently, the framework of the continuous time model---the point process---
has been gradually applied to model financial data at the transaction
level \citep{Hasbrouck:1991,Engle:1998,Bowsher:2007,Hawkes:1971,Hawkes:1971-2,Hawkes:1974,Hawkes:2018,Filimonov:2012,Filimonov:2015,Wheatley:2019}. 

A point or counting process is characterized by its "intensity" function, which represents 
the conditional probability of an event occurring in the immediate future.
The Hawkes process, introduced by A. G. Hawkes \citep{Hawkes:1971,Hawkes:1971-2,Hawkes:1974}, is a type of point process 
originally developed to model the occurrence of 
seismic events. Owing to its simplicity and flexibility, this model is being increasingly used
in high-frequency finance \citep{Bacry:2015,Kirchner:2017,Blanc:2017,Errais:2010}.
It can easily capture the interactions between different types of events, incorporate the influence of intensive factors through marks, and accommodate non-stationarities.

The non-stationarity of the Hawkes process arise from the reproduction of an event 
induced by the interaction terms described by the kernel function that estimate the effect of the past events based on the elapsed time from them.
The intensity of the process is determined by summing the effects of all past events through the kernel function. 
Additionally, by introducing a real number "mark" that represents the strength of event effects, both the marks and the kernel values can be considered in the intensity estimation.
If the average number of events triggered by one event, known as the branching ratio, exceeds one, the system becomes non-stationary.
If the branching ratio approaches one from below, the system becomes critical, leading to power-law scaling in
the PDF of the intensities in the intermediate asymptotic region.
In the sub-critical case, the distribution decays exponentially beyond the characteristic intensity, 
which diverges as the branching ratio approaches one \citep{Kanazawa:2020,Kanazawa:2020-2}.

The discrete time self-exciting negative binomial distribution (DT-SE-NBD) was proposed for the modeling of time-series data related to defaults \citep{Hisakado:2020}.
This model incorporates a parameter $\omega$ that captures the 
correlation of events within the same term, and the process converges to a
discrete-time self-exciting Poisson process in the limit $\omega \to 0$.
In the context of time-series default data, defaults are typically recorded on a yearly or quarterly basis, and defaults within the same period exhibit high correlation. Consequently, the probability distribution of the number of defaults displays overdispersion, where the variance is significantly larger than the mean value. To accurately capture this overdispersion, the negative binomial distribution (NBD), which has two parameters to control both the mean and variance, is more suitable than the Poisson distribution    \citep{Hisakado:2022,Hisakado:2022-1}. 
In the continuous time limit, the 
DT-SE-NBD process transforms into a point process with marks. As $\omega$ tends to zero,
the marked point process becomes a "pure" Hawkes process as 
 the value of the marks equals one. 
The PDF of the intensities also exhibits power-law scaling and exponentially decays beyond the characteristic intensity 
at the critical point and in the subcritical case, respectively.

This study extends the results of the
SE-NBD process with a single exponential 
kernel to the case where the memory kernel is the sum of an arbitrary 
finite number of exponentials. 
Empirically, power-law memory 
kernels have been observed in studies on financial transaction data\citep{Bacry:2015}.
The theoretical analysis of Hawkes processes with power-law memory kernels requires decomposing the power-law memory kernel into a superposition of exponential kernels with weights following the inverse gamma distribution.
By considering the case of multiple exponential kernels,
we can estimate the power-law exponent of the intensity distribution for the marked Hawkes
process with a power-law memory kernel.
Additionally, we validate the theoretical predictions for the power-law behavior of the intensities through numerical simulations.
The remaining sections of this paper are organized as follows: In Section \ref{sec:model}, 
we introduce an SE-NBD process with an arbitrary finite number of exponentials
and review related results on the process as well as on the Hawkes process.
Section \ref{sec:solution} focuses on studying the multivariate stochastic differential equation for the process, providing the solution to the two-variate master equations with a double exponential memory kernel, and deriving the PDF of the intensities.
Furthermore, we discuss the PDF of the intensities for the case where the memory kernel is a sum of finite $K$ exponential functions. In Section \ref{sec:Numerical}, we verify the theoretical predictions for the power-law exponents of the PDF of the intensities at the critical point through numerical simulations. 
Finally, we present our conclusions in Section \ref{sec:conclusion}.

\section{Model}
\label{sec:model}
We consider a DT-SE-NBD process $\{X_t\},t=1,\cdots$.
The variable $X_t \in \{0,1,\cdots\}$ represents the size of the event at time $t$.
In the context of modeling of time-series default data, $X_t$ specifically represents the number of 
defaults that occur in the $t$-th period\citep{Hisakado:2022,Hisakado:2022-1}.
$X_t$ obeys NBD for the condition $\hat{\lambda}_t=\lambda_t$ as 
\begin{eqnarray}
X_{t+1}&\sim& \mbox{NBD}
\left(\alpha=\frac{\hat{\lambda}_{t}}{\omega },p=\frac{1}{\omega+1}\right),t\ge 0
\label{NBD}
\\
\hat{\lambda}_t&=&\nu_0+n\sum_{s=1}^{t}h_{t-s} X_s,t \ge 1\,\,,\,\,\hat{\lambda}_0=\nu_0 
\label{intensity}
\\
h_t&=&\frac{1}{n}
\sum_{k=1}^K n_k h^k_t =\frac{1}{n}\sum_{k=1}^{K}n_k (1-e^{-1/\tau_k})e^{-t/\tau_k},
n=\sum_{k=1}^{K} n_k \label{dmemory}. 
\end{eqnarray}
Here, $\alpha,p$ are the parameters of the NBD. 
$\alpha > 0$ 
represents the number of successes until the experiment is terminated, 
and $p\in(0,1]$ is the probability of success in each individual experiment.
The sequence $\{h_t\}_{t=0,1,\cdots}$ refers to the discount factor (memory kernel), which is assumed to be normal; specifically, $\sum_{t=0}^{\infty}h_t=1$. 
The prefactor $(1-e^{-/\tau_k})$ is introduced to ensure the normalization 
of the $k$-th exponential term $h^k_t=(1-e^{-1/\tau_k})e^{-t/\tau_k}$, where $\tau_k$ represents the memory length associated with the $k$-th term. Consequently, $\sum_{t=0}^{\infty}(1-e^{-1/\tau_k})e^{-t/\tau_k}=1$.
The set of coefficients $\{n_k\}_{k=1,\cdots,K}$ 
quantifies the contribution of the $k$-th 
exponential term to the memory kernel $h_t,t=0,1\cdots$, with
$n$ representing their sum.

\begin{figure}[htbp]
\begin{center}
\includegraphics[width=120mm]{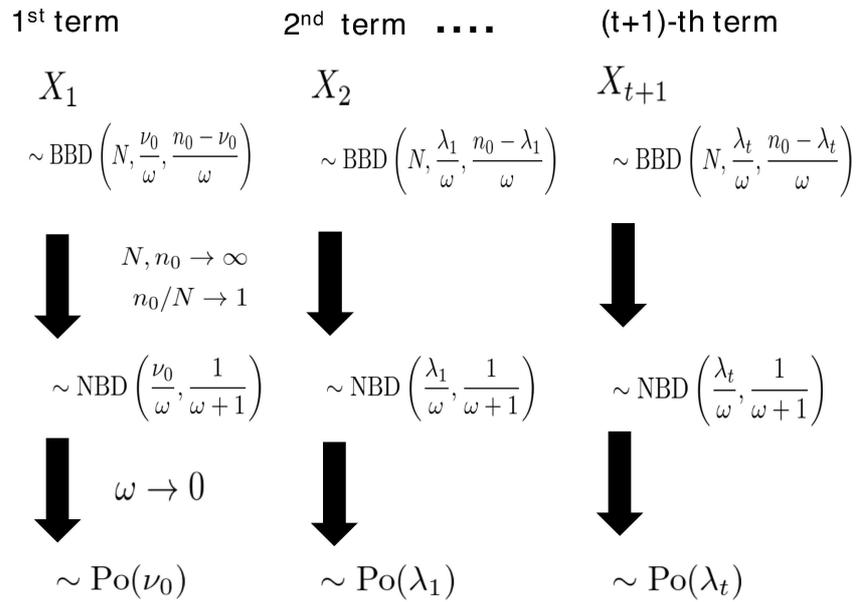}
\caption{Transformation of multi-term P\`{o}lya urn process (beta-binomial) to discrete-time Negative binomial and  Poisson process.}
\label{fig:hawkes}
\end{center}
\end{figure}

By considering a double-scaling limit, the DT-SE-NBD process can be constructed as an extension of the multi-term P\'{o}lya urn process with a memory kernel as follows:
\[
X_{t}\sim \lim_{\substack{N,n_{0}\to\infty \\ n_{0}/N=1}}\mbox{BBD}\left(N,\alpha=\frac{\hat{\lambda}_t}{\omega},\beta=\frac{n_0-\hat{\lambda}_t}{\omega}\right)=\mbox{NBD}\left(\alpha=\frac{\hat{\lambda}_{t}}{\omega },p=\frac{1}{\omega+1}\right).
\]
Here, BBD represents the beta-binomial distribution, and the variables $N,\hat{\lambda}_t,n_{0}-\hat{\lambda}_t$ represent the number of trials, the number of red balls, and the number of blue balls in an urn in the $t$-th term, respectively. 
In each term, we sequentially remove a ball and 
return it along with $\omega$ additional
balls of the same color. 
Consequently, the total number of balls in the urn
increases by $\omega$. The process is repeated $N$ times within one term.
The term $1/(1+n_0/\omega)$ represents Pearson's correlation coefficient, , which reflects the correlation between the color choices of the balls in the process \citep{Hisakado:2006}.
The feedback term $h_t$ induces an intertemporal correlation of the
ball choices. Furthermore, as $\omega$ approaches zero, 
$X_{t+1}$ under the condition $\hat{\lambda}_t=\lambda_t$ follows a Poisson random variable 
with a mean of $\lambda_t$, thereby transforming the DT-SE-NBD process into a  DT-SE Poisson process. 
\[
X_{t+1}\sim \mbox{Po}(\hat{\lambda}_{t}),t\ge 0
\]
Fig.\ref{fig:hawkes} illustrates the transformation from the multi-term P\`{o}lya urn process
 to the DT-SE Negative binomial and DT-SE Poisson process.

We study the (unconditional) expected value of $X_{t+1}$.
As $X_{t+1}\sim \mbox{NBD}
\left(\frac{\lambda_t}{\omega},\frac{1}{\omega+1}\right)$ for $\hat{\lambda}_t=\lambda_t$, 
we have
\[
E[X_{t+1}]=E_{\lambda_t}[E[X_{t+1}|\hat{\lambda}_t=\lambda_t]]=E_{\lambda_t}[\hat{\lambda}_t].
\]
where $E_{\lambda_t}[\,\,\,\,]$ is the ensemble average of the stochastic process 
$\{X_{s}, s \leq t\}$
[i.e., Filtration $F_t$]. We are now interested in the steady state of the process and write the unconditional
expected value of $X_t$ as $\mu$. We have the following relationship:
\[
\mu\equiv \lim_{t\to\infty}E_{\lambda_t}[\hat{\lambda}_t].
\]
With $\lambda_t= \nu_0 +n\sum_{s=1}^{t}X_s h_{t-s}$,$E[X_s]=\mu$ and the normalization of the discount factor $h_t$, we have
\[
\mu = \nu_0 + n\mu.
\]
We then obtain,
\[
    \mu=\frac{\nu_0}{1-n}
\]
The phase transition between the steady and non-steady state occurs at
$n = n_c = 1$, which is the critical point. This is the same as that
in the Hawkes process. In this sense, $n$ 
corresponds to the branching ratio of the Hawkes process, and
the steady-state condition is $n<1$.

We proceed to the continuous time SE-NBD process.
We begin with the DT-SE-NBD process in \eqref{NBD}, \eqref{intensity}, and \eqref{dmemory}. The stochastic process $\{X_t\},t=1,\cdots$ is non-Markovian, so we focus on 
the time evolution of $\hat{z}_t\equiv 
\hat{\lambda}_t-\nu_0$.
We decompose $\hat{z}_t$ as follows:
\begin{eqnarray}
\hat{z}_t&\equiv& \hat{\lambda}_t-\nu_0=\sum_{k}\hat{z}^k_{t}  \nonumber \\
\hat{z}^k_t&\equiv& n_k\sum_{s=1}^{t}h^k_{t-s}X_s\,\, , \,\,\ \hat{z}_{0}^k=0 \color{black}.
\nonumber
\end{eqnarray}
$\hat{z}^k_t$ satisfies the following recursive equation:
\[
\hat{z}^k_{t+1}=e^{-1/\tau_k}\hat{z}^{k}_t+ n_k h^{k}_0 X_{t+1}.
\]
Here, we use the relation $\sum_{s=1}^{t+1}X_s h^k_{t+1-s}=X_{t+1}h^k_0+e^{-1/\tau_k}\sum_{s=1}^{t}X_s h^k_{t-s}$.
The stochastic difference equation for $\hat{z}^{k}_{t}$ is
\begin{equation}\label{MSDE}
\Delta \hat{z}^k_t\equiv \hat{z}^{k}_{t+1}-\hat{z}^k_{t}
=(e^{-1/\tau_k}-1)\hat{z}^k_t+n_k h^k_0 X_{t+1},k=1,\cdots,K.
\end{equation}
The DT-SE-NBD process is now represented by a multivariate 
stochastic difference equation.

Now, we consider the continuous time limit. We divide the unit time interval 
by infinitesimal time intervals with width $dt$. The decreasing factor 
$e^{-1/\tau_k}$ for unit time period
is replaced by $e^{-dt/\tau_k}\simeq 1-dt/\tau_k+o(dt/\tau_k)$
 for $dt$. In addition, we adopt the 
notation $h_{k}(t)=\frac{1}{\tau_k}e^{-t/\tau_k}$. Here, we replace the normalization factor
$(1-e^{-1/\tau_k})$ of $h^k_t$ with $1/\tau_k$ in $h_{k}(t)$
to ensure that  $\int_{0}^\infty
h_{k}(t)dt=\frac{1}{\tau_k}\int_{0}^{\infty}e^{-t/\tau_k}=1$. We 
introduce the continuous time memory kernel $h(t)$ as
\begin{equation}\label{memory}
n h(t)=\sum_{k=1}^{K}n_k h_k(t)=\sum_{k=1}^K n_k\frac{1}{\tau_k} e^{-t/\tau_k}.
\end{equation}
We denote the  continuous time 
limits of $\hat{z}^{k}_{t}$ and $\hat{\lambda}_t$ as
$\hat{z}_k(t)$ 
and $\hat{\lambda}(t)$, respectively.
The timing, size of events, and counting process
are denoted by $\hat{t}_i$, $m_i$ and $\hat{N}(t)$, respectively.
$\hat{z}_k(t)$  and $\hat{\lambda}(t)$ are defined as
\begin{eqnarray}
\hat{z}_{k}(t)&=&n_{k}\sum_{i=1}^{\hat{N}(t)}h_{k}(t-\hat{t}_i)\cdot m_i \nonumber \\
\hat{\lambda}(t)&=&\nu_0+\sum_{k=1}^{K}\hat{z}_{k}(t) \nonumber 
\end{eqnarray}

$X_{t}$ is the common noise for the time interval $[t,t+1)$ in the DT SE-NBD process. To take the continuous time limit, it is necessary to define the noise $d\xi(t)$ for the infinitesimal interval $[t,t+dt)$. There is freedom in defining $d\xi(t)$, and we adopt the following, which is based on the reproduction properties of NBD.
When $X_1\sim NBD(\alpha_1,p),X_2\sim \mbox{NBD}(\alpha_2,p)$,$X_1+X_2\sim \mbox{NBD}(\alpha_1+\alpha_2,p)$.
We define $d\xi^{NBD}_{(\alpha,p)}(t)$  for the interval $[t,t+dt)$
as
\[
d\xi^{NBD}_{(\alpha,p)}(t) \sim \mbox{NBD}(\alpha dt,p).
\]
If we consider the above-mentioned noise, by the 
reproduction property, the noise for the interval $[t,t+1)$ becomes
\[
\int_{t}^{t+1}d\xi^{NBD}_{(\alpha,p)}(s)\sim \mbox{NBD}(\alpha,p).
\]
In the self-exciting process, we replace $\alpha$ and $p$ by $\frac{\hat{\lambda}(t)}{\omega}$ and $\frac{1}{\omega+1}$, respectively.
The conditional expected value and variance of $d\xi^{NBD}_{(\alpha=\frac{\hat{\lambda}(t)}{\omega}dt,p=\frac{1}{\omega+1})}$ are
\begin{eqnarray}
&&E\left[d\xi^{NBD}_{(\frac{\hat{\lambda}(t)}{\omega}dt,\frac{1}{\omega+1})}\right|\left.\hat{\lambda}(t)=\lambda(t)\right]=\lambda(t)dt \nonumber \\
&&V\left[d\xi^{NBD}_{(\frac{\hat{\lambda}(t)}{\omega},\frac{1}{\omega+1})}\right|\left.\hat{\lambda}(t)=\lambda(t)\right]=\lambda(t)(\omega+1)dt \nonumber .
\end{eqnarray}

The conditional probability of occurrence of an event with size $m$ during the time interval
$[t,t+dt)$ is given as
\[
P\left(d\xi^{NBD}_{(\frac{\hat{\lambda}(t)}{\omega},\frac{1}{\omega+1})}=m\right|\left.\hat{\lambda}(t)=\lambda(t)\right)
=\left\{
\begin{array}{cc}
1-\frac{\lambda(t)}{\omega}\log (\omega+1)dt  & m=0 \\
\frac{1}{m}\left(\frac{\lambda(t)}{\omega}\right)
\left(\frac{\color{red}\omega\color{black}}{\omega+1}\right)^m dt & m\ge 1.
\end{array}
\right.
\]
In the limit $\omega\to 0$, the probabilities converge to
\[
\lim_{\omega\to 0}
P\left(d\xi^{NBD}_{(\frac{\hat{\lambda}(t)}{\omega},\frac{1}{\omega+1})}=m\right|\left.\hat{\lambda}(t)=\lambda(t)\right)
=\left\{
\begin{array}{cc}
1-\lambda(t) dt  & m=0 \\
\lambda(t)dt & m= 1 \\
o(dt) & m\ge 2,
\end{array}
\right.
\]
and the event size $m$ is restricted to one.

The NBD noise $\xi^{NBD}(t)$ is then written as
\[
\xi^{NBD}_{\left(\frac{\hat{\lambda}(t)}{\omega},\frac{1}{\omega+1}\right)}
(t)=\sum_{i=1}^{\hat{N}(t)}m_i \delta(t-\hat{t}_i).
\]
The probability mass
function (PMF) for event size
$m=1,\cdots$ is given by
\begin{equation}\label{PMF}
\rho(m)=\frac{P\left(d\xi^{NBD}_{(\frac{\hat{\lambda}(t)}{\omega},\frac{1}{\omega+1})}=m\right)}{P(d\xi^{NBD}>0)}=\frac{1}{m \log (\omega+1)}\left(\frac{\color{red}\omega\color{black}}{\omega+1}\right)^m.
\end{equation}
The intensity of the process, $P(d\xi^{NBD}>0)/dt$,  is then given as
\begin{equation}
\hat{\nu}(t)\equiv \frac{\log (\omega+1)}{\omega}\hat{\lambda}(t)\,\, , \,\, \hat{\lambda}(t)=\nu_0+\sum_{i=1}^{\hat{N}(t)}\sum_{k=1}^{K}n_{k} h_{k}(t-\hat{t}_i)\cdot m_i \label{eq:intensity}.
\end{equation}
This process is known as the marked Hawkes process. Each event has a distinct 
influential power that is encoded in the marks $\{m_i\}_{i=1,\cdots,\hat{N}(t)}$. 
$\{m_{i}\}_{i=1,\cdots,\hat{N}(t)}$ are IID
random numbers and obey the PMF of \eqref{PMF}.
When the parameter $\omega$ tends to zero ($\omega \to 0$), the PMF $\rho(m)$ approaches the delta function $\delta(m-1)$, and the marked Hawkes process reduces to the “pure” Hawkes process.
In the following analysis, we will focus on studying the PDF of the intensities in the marked Hawkes process.

\section{Solution}
\label{sec:solution}
The multivariate difference equation \eqref{MSDE}
can be transformed into a multivariate SDE as follows:
\begin{eqnarray}
\hat{z}_{k}(t)&=&-\frac{1}{\tau_k}\hat{z}_k(t)dt+\frac{n_k}{\tau_k}d\xi^{NBD}_{(\frac{\hat{\lambda}(t)}{\omega},\frac{1}{\omega+1})},k=1,\cdots,K.\label{SDE2} \\
\hat{\lambda}(t)&=&\nu_0+\sum_{k=1}^{K}\hat{z}_{k}(t) \label{nu_t}\\
\xi^{NBD}_{\left(\frac{\hat{\lambda}(t)}{\omega},\frac{1}{\omega+1}\right)}
(t)&=&\sum_{i=1}^{\hat{N}(t)}m_i \delta(t-\hat{t}_i) \label{noise2}.
\end{eqnarray}
Note that the same state-dependent NBD noise
$\xi^{\mathrm{NBD}}_{(\frac{\hat{\lambda}}{\omega},\frac{1}{\omega+1})}$
affects every component of the multivariate SDE $\{\hat{z}_k\}_{k=1,...K}$.
In other words, each shock event simultaneously affects the trajectories for all excess intensities $\{\hat{z}_k\}_{k=1,...K}$.

The formal solution of the SDE is
\[
\hat{z}_k(t)=\frac{n_k}{\tau_k}\int^{t}_0dt'e^{-(t-t')/\tau_k}\xi^{\mathrm{NBD}}_{(\frac{\hat{\lambda}}{\omega},\frac{1}{\omega+1})}(t').
\]
The SDE\eqref{SDE2} can be interpreted as
\[
\hat{z}_k(t+dt)-\hat{z}_k(t)=
\left\{ 
\begin{alignedat}{2}   
    -\frac{\hat{z}_k(t)}{\tau_k}dt \qquad &\mathrm{Prob.}=1-\frac{\hat{\lambda}(t)}{\omega}\log(\omega+1)dt  \\
    \frac{n_km}{\tau_k} \hspace{ 25pt } &\mathrm{Prob.}=\frac{1}{m}\left(\frac{\lambda(t)}{\omega}\right)\left(\frac{\color{red}\omega\color{black}}{\omega+1}\right)^mdt,m\geq1  .
\end{alignedat} 
\right.
\]

We adopted the same procedure to derive the master equation for the PDF of $\hat{z}_k$ in Ref.\citep{Kanazawa:2020,Kanazawa:2020-2}. 
As the SDEs for $\boldsymbol{z}:=(z_1,...,z_K)$ are standard Markovian stochastic processes, we obtain the corresponding master equation
\begin{align}\label{master}
\frac{\partial}{\partial t}P_{t}&(\boldsymbol{z})=\sum^K_{k=1}\frac{\partial}{\partial z_k}\frac{z_k}{\tau_k}P_t(\boldsymbol{z})+\sum_{m=1}^{\infty}\frac{1}{m\omega}\left(\frac{\omega }{\omega+1}\right)^m \notag\\
&\times \left\{\left[\nu_0+\sum ^K_{k=1}\left(z_k-\frac{n_km}{\tau_k}\right)\right]P_t(\boldsymbol{z-h})-\left[\nu_0+\sum^K_{k=1}z_k\right]P_t(\boldsymbol{z})\right\}
\end{align}
The jump-size vector is given by $\boldsymbol{h}:=(n_1m/\tau_1,...,n_Km/\tau_K)$. 

The master equation \eqref{master} takes a simplified form under the Laplace representation:
\begin{equation}\label{laplace:z}
    \tilde{P}_t(\boldsymbol{s}):=\mathcal{L}_{\color{red}K\color{black}}[P_t(\boldsymbol{z});\boldsymbol{s}],
\end{equation}
where the Laplace transformation in the $K$-dimensional space is defined as
\begin{equation}
    \mathcal{L}_K[P_t(\boldsymbol{z});\boldsymbol{s}]:=\int^{\infty}_0 d\boldsymbol{z} e^{-\boldsymbol{s\cdot z}}P_t(\boldsymbol{z})
\end{equation}
with the volume element $d\boldsymbol{z}:=\prod^K_{k=1}dz_k$. The wave vector $\boldsymbol{s}:=(s_1,...,s_K)$ is the conjugate of the excess intensity vector 
$\boldsymbol{z}:=(z_1,...,z_K)$.

The Laplace representation of the master equation \eqref{master} is given by
\begin{equation}\label{laplace:master}
    \frac{\partial \tilde{P}_t(\boldsymbol{s})}{\partial t} = -\sum^K_{k=1}\frac{s_k}{\tau_k}\frac{\partial \tilde{P}_t(\boldsymbol{s})}{\partial s_k}+
    \sum^{\infty}_{m=1}\frac{1}{m\omega}\left( \frac{\omega} {\omega+1} \right)^m\left(\nu_0-\sum^K_{k=1}\frac{\partial}{\partial s_k}\right)(e^{-\boldsymbol{h\cdot s}}-1)\tilde{P}_t(\boldsymbol{s})
\end{equation}
Then, the Laplace representation \eqref{laplace:z} of $P_t(\boldsymbol{z})$,  which is the solution of \eqref{laplace:master}, enables us to obtain the (one-dimensional) Laplace representation $\tilde{Q}_t(s)$ of the intensity $\mathrm{PDF}\ P_t(\nu)$ according to 
\begin{equation}\label{laplace:nu}
    \tilde{Q}_t(s):=\mathcal{L}[P_t(\nu);s]=e^{-\nu_0s}\tilde{P}_t(\boldsymbol{s}=(s,s,...,)).
\end{equation}

\subsection*{A. Single exponential kernel}
We now consider the case in which the memory function \eqref{memory}
consists of $K=1$ exponential functions\citep{Hisakado:2022}. 
The Laplace representation of the master equation \eqref{master} is given by
\[
    \frac{d \tilde{P}_t(s)}{d t} = -\frac{s}{\tau}\frac{d \tilde{P}_t(s)}{ds}+
    \sum^{\infty}_{m=1}\frac{1}{m\omega} \left( \frac{\omega} {\omega+1} \right)^m\left(\lambda_0-\frac{d}{d s} \right)\left(e^{-\frac{nk}{\tau}s}-1\right)\tilde{P}_t(s).
\]
The steady-state PDF of the intensity $\hat{\lambda}(t)$ is
\begin{equation}\label{power-law_exponent:1}
P_{SS}(\lambda)\propto \lambda^{-1+2\frac{\nu_0 \tau}{\omega +1}}e^{-\frac{2\tau\epsilon}{\omega +1}\lambda}\,\, ,\,\, \epsilon=1-n. 
\end{equation}
The power-law exponent of the PDF of the excess intensity is $1-\frac{2\nu_0\tau}{\omega+1}$ and depends on $\omega$. In the limit $\omega \rightarrow 0$, the result coincides with that in Ref.\citep{Kanazawa:2020,Kanazawa:2020-2}.
The power-law exponent increases with $\omega$ and converges to $1$ in the limit $\omega \rightarrow\infty$. In addition, the length scale beyond which the intensity shows
exponential decay for the off-critical case is $(\omega+1)/(2\tau\epsilon)$,
and it diverges in the limit $\epsilon=1-n \to 0$.
This is also an increasing function of $\omega$.

\subsection*{B. Double exponential kernel}
We consider the case where the memory
function \eqref{memory} consists of $K=2$ exponential functions.
To derive the solution for the Laplace representation of the master equation \eqref{laplace:master},
we apply the method of characteristics as in Ref.\citep{Kanazawa:2020-2}. One can find a brief explanation 
in Refs.\citep{gardiner} and \citep{Kanazawa:2020-2}. In appendix \ref{sec:characteristics}, we provide a brief review.

We start from the Lagrange--Charpit equations, which are given by
\begin{eqnarray}
    \frac{ds_1}{dl}&=&-\sum^{\infty}_{m=1}\frac{1}{m\omega} \left( \frac{\omega} {\omega+1} \right)^m\left(e^{-\boldsymbol{h\cdot s}}-1\right)-\frac{s_1}{\tau_1},\label{a}\\
    \frac{ds_2}{dl}&=&-\sum^{\infty}_{m=1}\frac{1}{m\omega} \left( \frac{\omega} {\omega+1} \right)^m\left(e^{-\boldsymbol{h\cdot s}}-1\right)-\frac{s_2}{\tau_2},\label{b}\\
    \frac{d}{dl}\log\tilde{P}_{ss}&=&-\sum^{\infty}_{m=1}\frac{1}{m\omega} \left( \frac{\omega} {\omega+1} \right)^m\nu_0\left(e^{-\boldsymbol{h\cdot s}}-1\right)\label{c}\\
\end{eqnarray}
and $l$ is the auxiliary “time” parameterizing the position on the characteristic curve. Let us develop the stability
analysis around $s = 0$ (i.e., for large $\lambda's$) for this pseudo-dynamical system.

\subsubsection*{a. Sub-critical case $n<1$}
\label{sub-critical}
Assuming $n:=n_1+n_2<1$, we first expand $e^{-\boldsymbol{h\cdot s}}$ to compute the sum of $m$:
\begin{equation}\label{Maclaurin}
    e^{-\left(\frac{n_1s_1}{\tau_1}+\frac{n_2s_2}{\tau_2}\right)m}\simeq1-\left(\frac{n_1s_1}{\tau_1}+\frac{n_2s_2}{\tau_2}\right)m+\frac{1}{2}\left(\frac{n_1s_1}{\tau_1}+\frac{n_2s_2}{\tau_2}\right)^2m^2+\cdots
\end{equation}
We obtain the linearized dynamics of the system \eqref{a},\eqref{b},\eqref{c} as follows:
\begin{equation}\label{eq:linearized}
    \frac{d\boldsymbol{s}}{dl}\simeq-\boldsymbol{H}\boldsymbol{s} ,\ 
    \frac{d}{dl}\log\tilde{P}_{ss}\simeq\nu_0\boldsymbol{K}\boldsymbol{s}
\end{equation}
with
\[
    \boldsymbol{H}:=
    \begin{pmatrix}
    \frac{1-n_1}{\tau_1} & \frac{-n_2}{\tau_2}\\
    \frac{-n_1}{\tau_1} & \frac{1-n_2}{\tau_2}
    \end{pmatrix}, \
    \boldsymbol{K}:=\left(\frac{n_1}{\tau_1},\frac{n_2}{\tau_2}\right).
\]
We introduce the eigenvalues $\beta_1,\beta_2$ and eigenvectors $\boldsymbol{e}_1,\boldsymbol{e}_2$ of $\boldsymbol{H}$ such that
\[
    \boldsymbol{P}:=(e_1,e_2),\ 
    \boldsymbol{P}^{-1}\boldsymbol{H}\boldsymbol{P}=
    \begin{pmatrix}
    \beta_1 & 0 \\
    0 & \beta_2
    \end{pmatrix}.
\]
The matrix $\boldsymbol{H}$ is the same  with the results of Ref.\citep{Kanazawa:2020,Kanazawa:2020-2} and  all eigenvalues are real. We denote
 them as $\beta_1\geq \beta_2$. The determinant of $\boldsymbol{H}$ is given by
\[
    \mathrm{det}\boldsymbol{H}=\frac{1-n}{\tau_1\tau_2}.
\]
This implies that the zero eigenvalue $\beta_1=0$ appears at the critical point $n=1$. Below the critical point $n<1$, all
eigenvalues are positive ($\beta_1,\beta_2>0$). For $n<1$, the dynamics can be rewritten as
\[
    \frac{d}{dl}\boldsymbol{P}^{-1}\boldsymbol{s}=-
    \begin{pmatrix}
    \beta_1 & 0 \\
    0 & \beta_2
    \end{pmatrix}
    \boldsymbol{P}^{-1}\boldsymbol{s}
    \Longrightarrow
    \boldsymbol{s}(l)=
    \boldsymbol{P}
    \begin{pmatrix}
    e^{-\beta_1(l-l_0)} \\
    e^{-\beta_2(l-l_0)}/C_1
    \end{pmatrix}.
\]
Here, $l_0$ and $C_1$ are integration constants.
We can assume $l_0=0$ as the initial point of the characteristic curve without loss of generality. Integrating the second equation in \eqref{eq:linearized}, we obtain
\begin{eqnarray}
    \log \tilde{P}_{ss}&=&\nu_0\boldsymbol{K}\int \boldsymbol{s}(l)dl+C_2=-\nu_0\boldsymbol{KP}
    \begin{pmatrix}
    1/\beta_1 & 0 \\
    0 & 1/\beta_2
    \end{pmatrix} \boldsymbol{P}^{-1}\boldsymbol{s}+C_2\nonumber \\
    &=&-\nu\boldsymbol{KH}^{-1}\boldsymbol{s}+C_2 \nonumber.
\end{eqnarray}
The general solution is given by
\begin{equation}\label{general_solution}
    \mathcal{H}(C_1)=C_2
\end{equation}
with function $\mathcal{H}$ determined by the initial condition of the characteristic curve. Let us introduce
\[
    \bar{\boldsymbol{s}}:=\boldsymbol{P}^{-1}\boldsymbol{s}=
    \begin{pmatrix}
    \bar{s}_1 \\
    \bar{s}_2
    \end{pmatrix}
    \Longrightarrow C_1=(\bar{s}_1)^{\beta_2/\beta_1}(\bar{s}_2)^{-1}.
\]
This implies that the solution is given in the following form:
\[
\log \tilde{P}_{ss}(\boldsymbol{s})=-\nu\boldsymbol{KH}^{-1}\boldsymbol{s}+\mathcal{H}((\bar{s}_1)^{\beta_2/\beta_1}(\bar{s}_2)^{-1}).
\]
Owing to the renormalization of the PDF, the relation 
\[
    \lim_{\boldsymbol{s\rightarrow0}} \log \tilde{P}_{ss}(\boldsymbol{s})=0
\]
must hold for any path in the ($s_1,s_2$) space ending at the origin ($\mathrm{limit}\  \boldsymbol{s\rightarrow0}$). Let us consider a specific limit such that $\bar{s}_2=x^{-1}(\bar{s}_1)^{\beta_2/\beta_1}$ and $\bar{s}_1\rightarrow 0$ for an arbitrary positive $x$.
\[
    \lim_{\bar{s}_1\rightarrow 0} \log\tilde{P}_{ss}(\boldsymbol{s})=\mathcal{H}(x)
\]
Because the left-hand side is zero for any $x$, the function $\mathcal{H}(\cdot)$ must be exactly zero. Thus, this leads to
\[
    \log\tilde{P}_{ss}(\boldsymbol{s})=-\nu\boldsymbol{KH}^{-1}\boldsymbol{s}.
\]
By substituting $\boldsymbol{s}=(s_1=s,s_2=s)$, we obtain
\[
    \log\tilde{Q}_{ss}(\boldsymbol{s})=-\nu_0s+\log\tilde{P}_{t=\infty}(s(1,1))
    \simeq-\frac{\nu_0}{1-n}s.
\]
This is consistent with the asymptotic mean intensity in the steady state \citep{Kanazawa:2020,Kanazawa:2020-2}.

\subsubsection*{b. Critical case $n=1$}
\label{critical}
In this case, the eigenvalues and eigenvectors of $\boldsymbol{H}$ are given by
\[
    \beta_1=0,\ \beta_2=\frac{n_1\tau_1+n_2\tau_2}{\tau_1\tau_2}, 
    \ \boldsymbol{e}_1=
    \begin{pmatrix}
    \tau_1 \\
    \tau_2
    \end{pmatrix},\ 
    \boldsymbol{e}_2=
    \begin{pmatrix}
    -n_2 \\
    n_1\end{pmatrix}.
\]
This means that the eigenvalue matrix and its inverse matrix are given by
\[
    \boldsymbol{P}=
    \begin{pmatrix}
    \tau_1 & -n_2\\
    \tau_2 & n_1
    \end{pmatrix},\ 
    \boldsymbol{P}^{-1}=\frac{1}{\alpha}
    \begin{pmatrix}
    n_1 & n_2 \\
    -\tau_2 & \tau_1
    \end{pmatrix},\ 
    \alpha := \det\boldsymbol{P}=\tau_1n_1+\tau_2n_2.
\]
Here, let us introduce
\begin{equation}\label{alpha}
    \boldsymbol{X}=\begin{pmatrix}
    X\\
    Y
    \end{pmatrix}=\boldsymbol{P}^{-1}\boldsymbol{s},
    \Longleftrightarrow X=\frac{n_1s_1+n_2s_2}{\alpha},\ 
    Y=\frac{-\tau_2s_1+\tau_1s_2}{\alpha}.
\end{equation}
We then obtain
\[
    \frac{dX}{dl}=0,\frac{dY}{dl}=-\beta_2Y
\]
at the leading linear order in the expansions of the powers of $X$ and $Y$. Because the first linear term is zero in the dynamics of $X$, corresponding to a transcritical bifurcation for the Lagrange--Charpit equations (16), we need to consider
the second-order term in $X$, namely,
\[
    e^{-\left(\frac{n_1}{\tau_1}s_1+\frac{n_2}{\tau_2}s_2\right)m}\simeq1-Xm+\frac{X^2m^2}{2}+n_1n_2\left(\frac{1}{\tau_1}-\frac{1}{\tau_2}\right)mY+\mathcal{O}(XY,X^2Y,Y^2)
\]
where we dropped the terms of the order $Y^2$, $XY$, and $X^2Y$. We then obtain the dynamic equations at the transcritical bifurcation to the leading order:
\[
    \frac{dY}{dl}\simeq-\beta_2Y,\  \frac{dX}{dl}\simeq-\frac{\omega+1}{2\alpha}X^2
\]
whose solutions are given by
\[
    X(l)=\frac{2\alpha}{\omega+1}\frac{1}{l-l_0}, \ Y(l)=C_1e^{-\beta_2(l-l_0)}
\]
with constants of integration $l_0$ and $C_1$. We can assume $l_0=0$ is the initial point on the characteristic curve.
Remarkably, only the contribution along the $X$ axis is dominant for the large $l$ limit (i.e., $|X|\gg|Y|$ for $l \rightarrow \infty$), which corresponds to the asymptotic $\boldsymbol{s}\rightarrow 0$. We then obtain
\begin{eqnarray}
\log\tilde{P}_{ss}&\simeq&\nu_0\int dl\left(\frac{n_1s_1(l)}{\tau_1}+\frac{n_2s_2(l)}{\tau_2}\right) \nonumber \\
    &\simeq& -\frac{2\nu_0\alpha}{\omega+1}\log X+\frac{\nu_0n_1n_2}{\beta_2}\left(\frac{1}{\tau_1}-\frac{1}{\tau_2}\right)Y+C_2 \nonumber 
\end{eqnarray}
with integration constant $C_2$. 
$C_2$ is a divergent constant because it has to compensate the diverging
logarithm $\log X$ to ensure that $\log P_{ss}(s=0) = 0$. This divergent
constant appears as a result of neglecting the ultraviolet (UV)
cutoff for small $s$ (which corresponds to neglecting the exponential tail 
of the PDF of intensities) \citep{Kanazawa:2020-2}.

Therefore, we obtain the steady solution
\[
    \log\tilde{P}_{ss}(\boldsymbol{s})=-\frac{2\nu_0\alpha}{\omega+1}\log X+\frac{\nu_0n_1n_2}{\beta_2}\left(\frac{1}{\tau_1}-\frac{1}{\tau_2}\right)Y
\]
for small $X$ and $Y$ by ignoring the UV cutoff and constant contribution. This recovers the power-law formula of the intermediate asymptotic of the PDF of the Hawkes intensities:
\begin{eqnarray}
\log\tilde{Q}_{ss}(s):&=&-\nu_0s+\log\tilde{P}_{ss}(s,s)\simeq-\frac{2\nu_0\alpha}{\omega+1}\log s\quad (s\sim0)\nonumber \\
    &\Longleftrightarrow& P(\lambda)\sim\lambda^{-1+\frac{2\nu_0\alpha}{\omega+1}} \quad (\lambda\rightarrow \infty) \label{power-law_exponent:2},
\end{eqnarray}
where $\alpha = \tau_1n_1+\tau_2n_2$, as defined in \eqref{alpha}. The power-law exponent of the PDF of the excess intensity is $1-\frac{2\nu_0\alpha}{\omega+1}$ and depends on $\omega$.
In the limit $\omega \rightarrow0$, the result coincides with that in Ref.\citep{Kanazawa:2020,Kanazawa:2020-2}. 

\subsection*{C. Discrete superposition of exponential kernels}
We now study the case in which the memory kernel is the sum of
$K$ exponentials for an arbitrary finite number $K$. We obtain the steady-state PDF of the intensity $\hat{\lambda}(t)$ as
\begin{equation}
  \label{power-law_exponent:K}
    \log\tilde{Q}_{ss}(s)\simeq-\frac{2\nu_0\alpha}{\omega+1}\log s\Longleftrightarrow P(\lambda)\sim\lambda^{-1+\frac{2\nu_0\alpha}{\omega+1}} ,\quad \mathrm{with} \ \alpha:=\sum^K_{k=1}n_k\tau_k.
\end{equation}
The power-law exponent of the PDF of the excess intensity is $1-\frac{2\nu_0\alpha}{\omega+1}$ and depends on $\omega$. In the limit $\omega \rightarrow 0$, the result coincides with that in Ref.\citep{Kanazawa:2020,Kanazawa:2020-2}. 

\subsection*{D. General case}
\label{General}
We study the case where the memory kernel is a continuous
superposition of exponential functions. We decompose the kernel as
\[
    h(t)=\frac{1}{n}\int^{\infty}_0 n(\tau)\frac{1}{\tau}e^{-t/\tau}d\tau\ ,\ n=\int^{\infty}_0 n(\tau)d\tau.
\]
This decomposition satisfies the normalization condition $\int^{\infty}_0 h(t)=1$. The function $n(\tau)$ quantifies the contribution of the exponential kernel $\frac{1}{\tau}e^{-t/\tau}$
to the branching ratio. The steady-state PDF of the intensity $\hat{\lambda}_t$ is
\begin{equation}\label{eq:alpha}
    \log\tilde{Q}_{ss}(s)\simeq-\frac{2\nu_0\alpha}{\omega+1}\log s\Longleftrightarrow P(\lambda)\sim\lambda^{-1+\frac{2\nu_0\alpha}{\omega+1}} ,\quad \mathrm{with} \ \alpha:=\int^{\infty}_0d\tau n(\tau)\tau. 
\end{equation}
The power-law exponent of the PDF of the excess intensity is $1-\frac{2\nu_0\alpha}{\omega+1}$ and depends on $\omega$. In the limit $\omega \rightarrow0$, the result coincides with that in Ref.\citep{Kanazawa:2020,Kanazawa:2020-2}. 

We now consider the case in which the memory kernel exhibits power-law decay. 
We express $h(t)$ as the superposition of $e^{-rt}$ with the weight of the 
gamma distribution $f_{\gamma}(r)=r^{\gamma-1}e^{-r}/\Gamma(\gamma)$ as
\[
    h(t)=\gamma(1+t)^{-(\gamma+1)}=\int^{\infty}_0re^{-rt}\cdot f_{\gamma}(r)dr, \ \gamma>0.
\]
By the change of variable $r\rightarrow1/\tau$,  we rewrite $h(t)$ as the superposition of $e^{-t/\tau}$ with the weight of the inverse gamma distribution
$f'_{\gamma}(\tau)=\tau^{-\gamma-1}e^{-1/\tau}/\Gamma(\gamma)$ as
\begin{equation}
    h(t)=\int^{\infty}_0\frac{1}{\tau}e^{-t/\tau}\cdot f_{\gamma}(1/\tau)\cdot\frac{1}{\tau^2}d\tau
        =\int^{\infty}_0n(\tau)\frac{1}{\tau}e^{-t/\tau}d\tau\ , n(\tau)=f'_{\gamma}(\tau) \label{eq:power}.
\end{equation}
$\alpha$ of Eq.\eqref{eq:alpha} is obtained as the expected value of the inverse 
gamma distribution,
\begin{equation}
\alpha = \int^{\infty}_0\tau f'_{\gamma}(\tau)d\tau
= \frac{1}{\gamma-1}, \ \gamma>1 \label{eq:power2}.
\end{equation}
The power-law exponent is $1-\frac{2\nu_0}{(\omega+1)(\gamma-1)}$.

\section{Numerical verification}
\label{sec:Numerical}
We conducted numerical simulations to study the steady-state PDF of the intensity in the marked Hawkes process. In particular, our goal was to verify the theoretical predictions of Eqs. \eqref{power-law_exponent:1}, \eqref{power-law_exponent:2}, \eqref{power-law_exponent:K}, and \eqref{eq:power2}. We employed the time-rescaling theorem, which enables us to sample data more efficiently than the rejection method. The time-rescaling theorem states that any point process with an integrable conditional intensity function can be transformed into a Poisson process with a unit rate \citep{Brown:2002}. That is, it is possible to convert a marked Hawkes process into a marked point process with a constant intensity of 1 by performing time-dependent rescaling. A point process with intensity 1 can easily generate the time of the next event occurrence because the event interval follows an exponential distribution.
The time-rescaling theorem introduces the following time transformations:
\begin{equation}\label{eq:time-rescaling}
    \Lambda(t) \equiv \int^{t}_0\hat{\nu}(s)ds.
\end{equation}

When event $\boldsymbol{t}_n=\{t_1,t_2,...,t_n\}$ follows a marked Hawkes process with intensity function $\hat{\nu}(t)$ in the observation interval $[0,T]$, event $\boldsymbol{t}'=\{\Lambda(t_1),\Lambda(t_2),...,\Lambda(t_n)\}$ obtained by the time transformation \eqref{eq:time-rescaling} follows a marked point process with intensity $1$ in the observation interval $[0,\Lambda(T)]$. The event interval $\Lambda(t_i)-\Lambda(t_{i-1})$ is given by
\[
\Lambda' \equiv \Lambda(t_i)-\Lambda(t_{i-1})=\int^{t_i}_{t_{i-1}}\hat{\nu}(s)ds \sim \mbox{Ex}(1).
\]
Hence, to numerically realize the marked Hawkes process, it is sufficient to generate a random number $\Lambda'$ following an exponential distribution with mean 1 and successively find $t_i$ that satisfies the above relation.

We rewrite the relation more concretely. For $t_{i-1}\le s \le t_{i}$, $z_{k}(s)$
is written as
\[
z_{k}(s)=n_{k}\sum_{j=1}^{i-1}\frac{1}{\tau_k}
e^{-(s-t_{j})/\tau_k}m_j=z_{k}(t_{i-1})e^{-(s-t_{i-1}/\tau_k}.
\]
The integral of $\nu(s)$ is 
\[
\int_{t_{i-1}}^{t_i}\nu(s)ds=\frac{\ln (\omega+1)}{\omega}\left(\nu_0(t_i-t_{i-1})+\sum_{k=1}^{K}z_{k}(t_{i-1})\int_{t_{i-1}}^{t_i}e^{-(s-t_{i-1})/\tau_k}ds\right).
\]
We have to solve the next relation,
\begin{eqnarray}
    \Lambda' &=& 
    \frac{\ln (\omega+1)}{\omega}\left(\nu_0(t_i-t_{i-1})+\sum_{k=1}^{K}z_{k}(t_{i-1})\tau_k(1-e^{-(t_{i}-t_{i-1})/\tau_{k}}\right) \nonumber \\
    &\Longrightarrow & \nu_0(t_{i}-t_{i-1})+\sum_{k=1}^{K}z_{k}(t_{i-1})\tau_k\left(1-\mbox{exp}\left(-\frac{t_i-t_{i-1}}{\tau_k}\right)\right)
    -\frac{\omega \Lambda'}{\ln(\omega+1)}=0  \label{eq:solve}
\end{eqnarray}
Note that the jump size of this process is determined according to the distribution of the mark $\rho(m)$ in Eq.\eqref{PMF} when the event occurrence time $t_{i}$ is determined.
The numerical procedures are given below.

\begin{figure}[htbp]
\begin{algorithm}[H]
    \caption{Sampling process of marked Hawkes process\\
    We denote $z_{k,i}=z_{k}(t_i)$.}
    \begin{algorithmic}[0]
    \STATE Step 1 $\hspace{5pt}$ $i\leftarrow 1$.
    \STATE Step 2 $\hspace{5pt}$ Get a random number $\Lambda'\sim \mbox{Ex}(1)$.
    \STATE Step 3 $\hspace{5pt}$ $t_1\leftarrow \frac{\Lambda'}{\nu_0}$.
    \STATE Step 4 $\hspace{5pt}$ Get a random number $m_i\sim \rho(m)$ in Eq.\eqref{PMF}.
    \STATE Step 5 $\hspace{5pt}$ $\hat{z}_{k,1}\leftarrow \frac{n_k m_1}{\tau_k},k=1,\cdots,K$.
    \STATE Step 6 $\hspace{5pt}$ The following steps are repeated.
    \STATE $\hspace{5pt}$ Step 6.1  $\hspace{5pt}$ $i\leftarrow i+1$.
    \STATE $\hspace{5pt}$ Step 6.2  $\hspace{5pt}$ Get a random number $\Lambda'\sim \mbox{Ex}(1)$.
    \STATE $\hspace{5pt}$ Step 6.3 $\hspace{5pt}$ Solve  Eq.\eqref{eq:solve} to get $t_i$.
    \STATE $\hspace{5pt}$ Step 6.4 $\hspace{5pt}$ Get a random number $m_i\sim \rho(m)$ in Eq.\eqref{PMF}.
    \STATE $\hspace{5pt}$ Step 6.5 $\hspace{5pt}$ $\hat{z}_{k,i}\leftarrow \hat{z}_{k,i-1}e^{-\frac{t_i-t_{i-1}}{\tau_k}}+\frac{n_k m_i}{\tau_k},k=1,\cdots,K$.
    \end{algorithmic}
\end{algorithm}
\end{figure}

We performed numerical simulations using Julia 1.7.3. The simulation time for $t\in [0,T=10^7]$ is approximately 0.5 h for background intensity $\nu_0=0.01$ and $K=1$.
The execution environment is as follows: OS: Ubuntu 20.045 LTS, Memory: 64 GB, CPU: 4.7 GHz.
The code is available on github\citep{github}. 
We sampled the point process $\{t_i\},i=1,...,N(T)$ with $T=10^7$. 
The common settings of the sampling processes are
$\nu_0\in\{0.01,0.2,1.0\},\omega\in\{0.01,1.0,10.0\},\tau=1.0$ and $n\in\{0.999,0.99,0.9\}$.

\subsubsection*{a. Single exponential kernel}
\begin{figure}[htbp]
\begin{center}
\begin{tabular}{ccc}
\includegraphics[width=5cm]{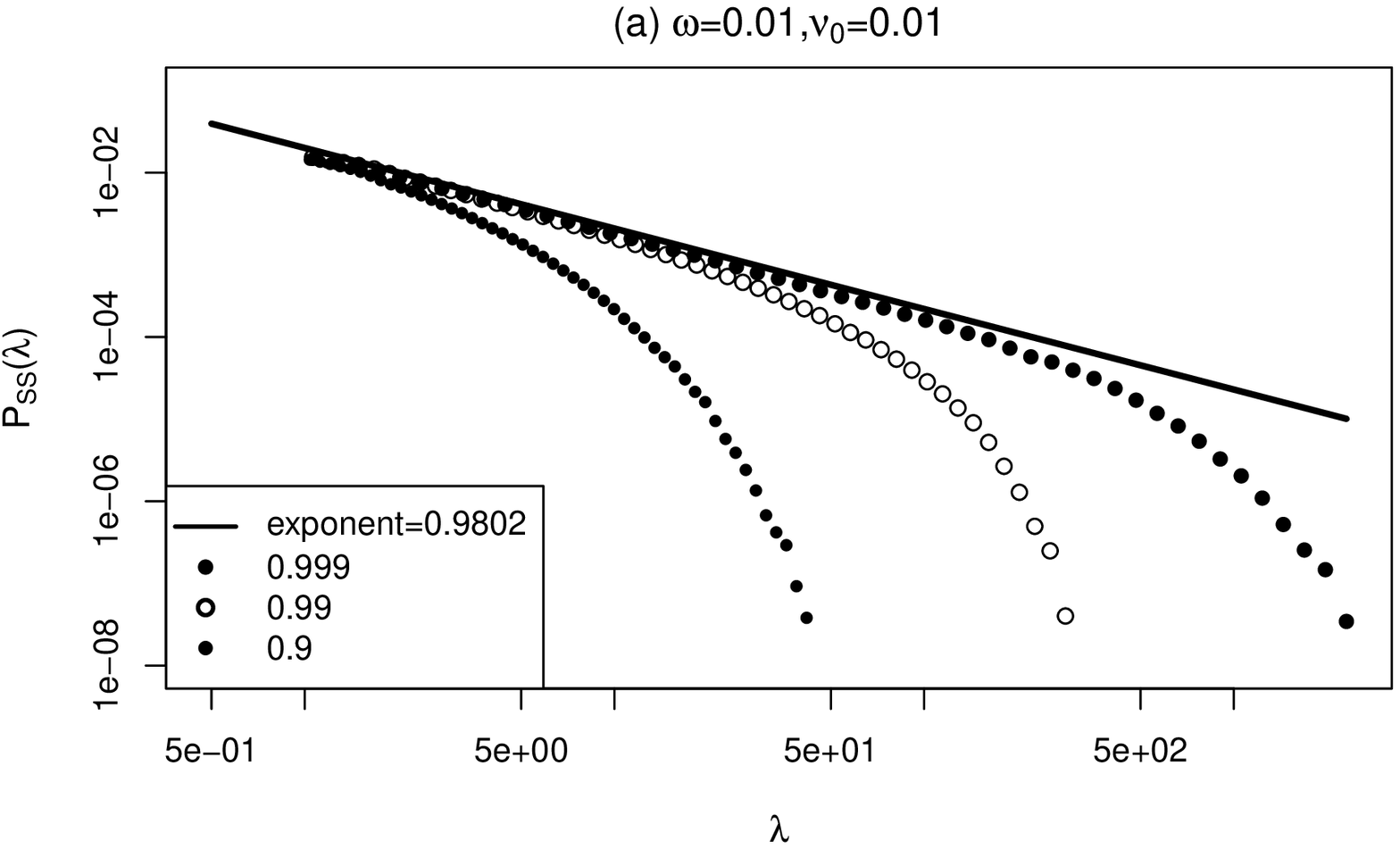}
&   
\includegraphics[width=5cm]{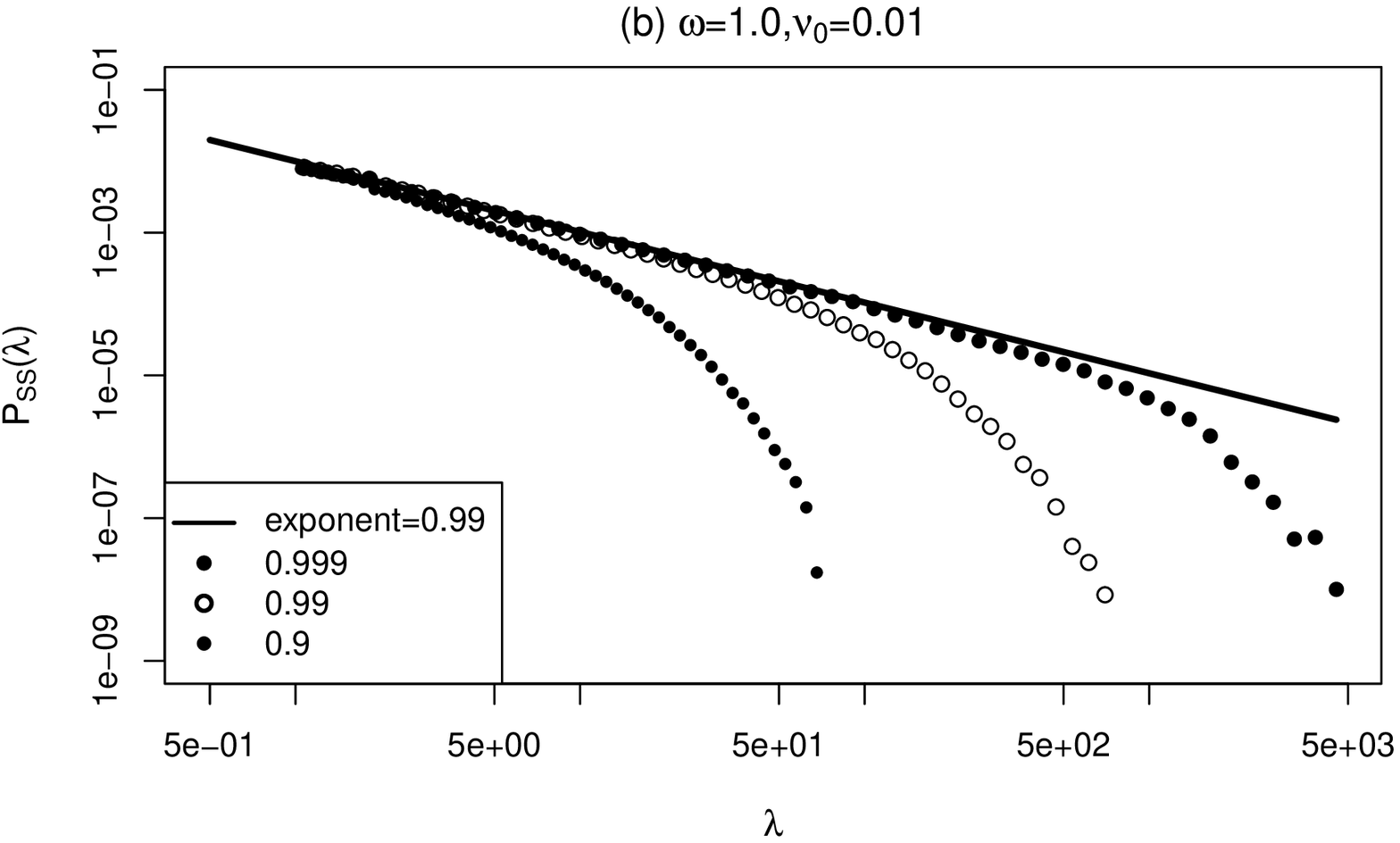}
&   
\includegraphics[width=5cm]{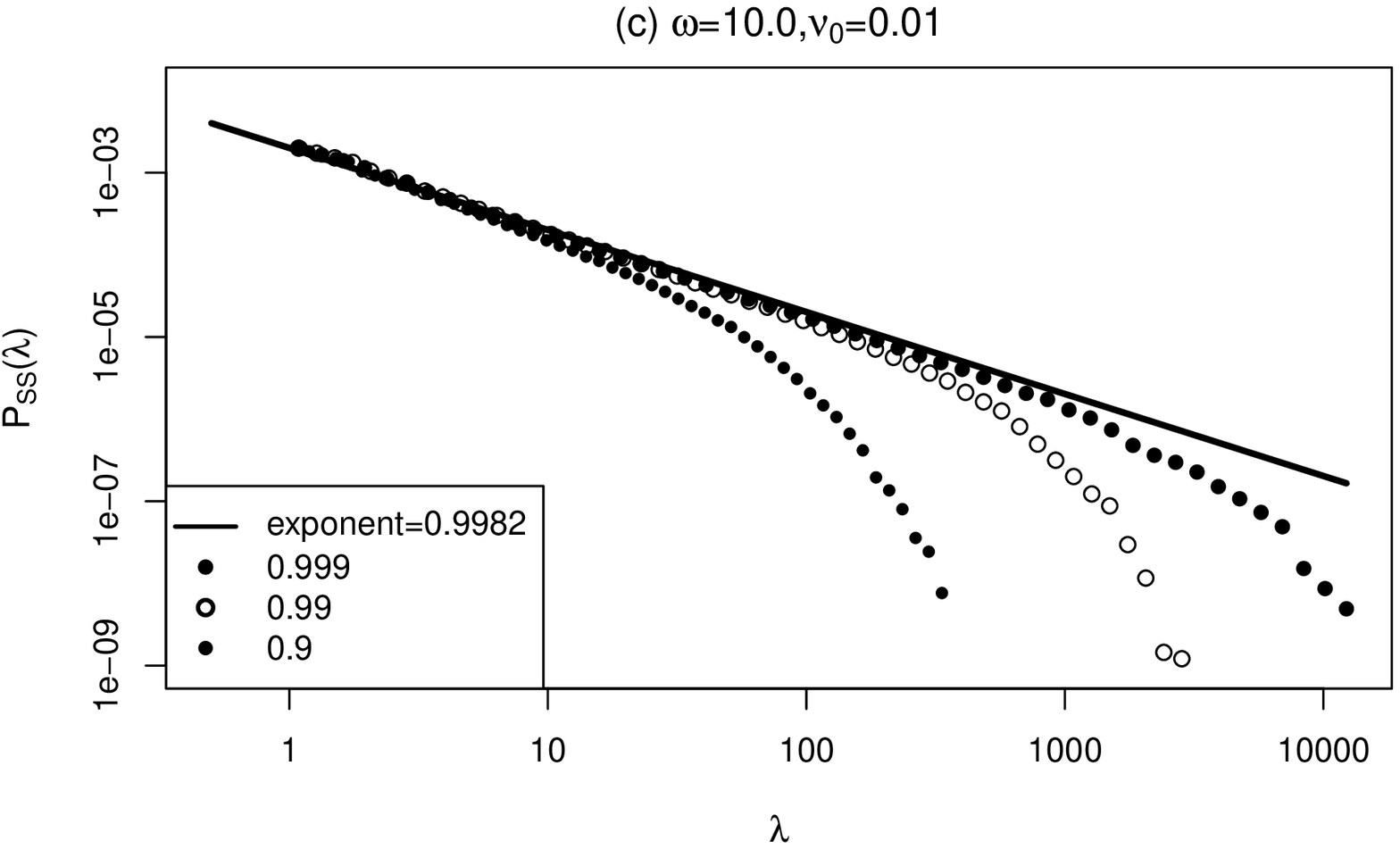}
\\
\includegraphics[width=5cm]{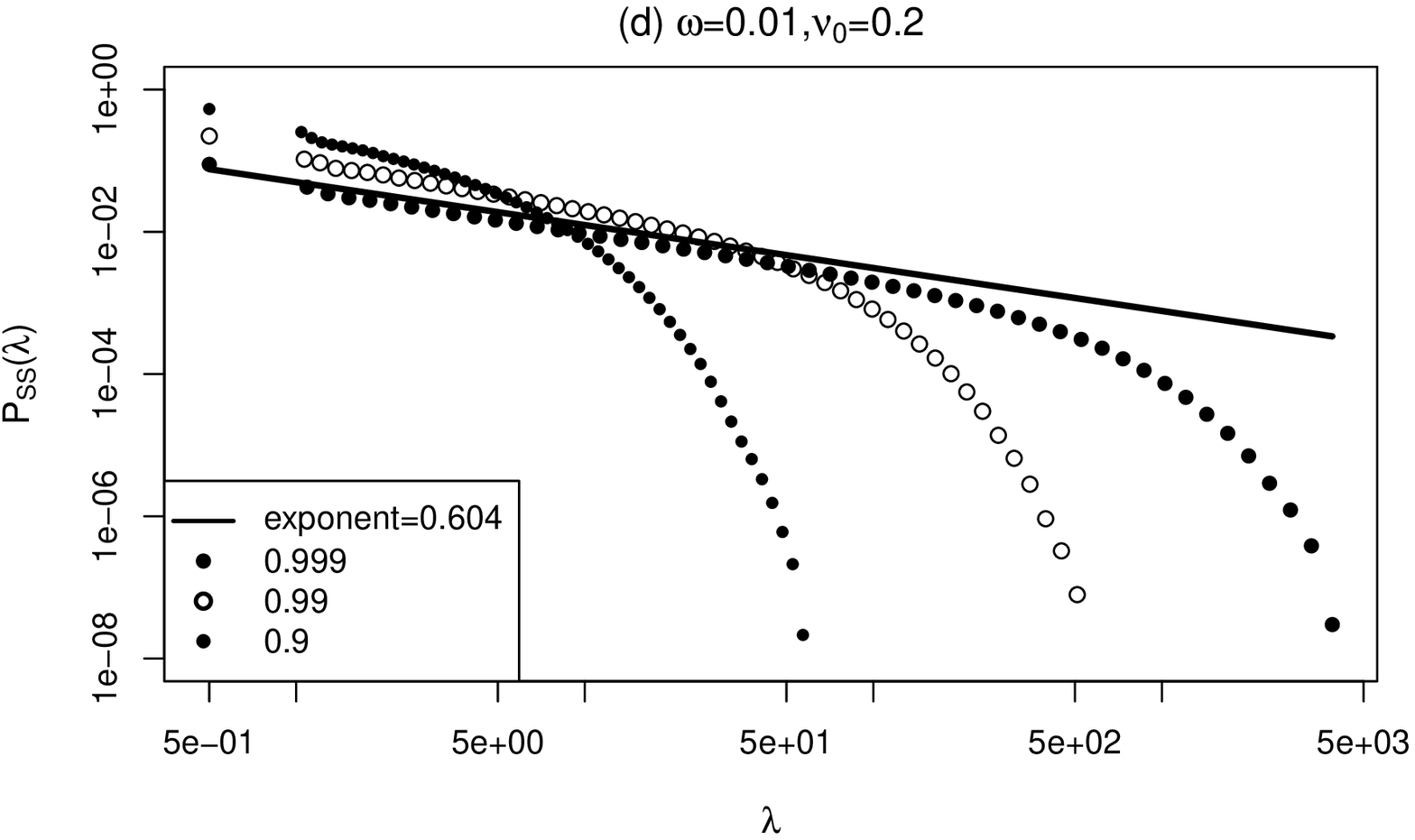}
&   
\includegraphics[width=5cm]{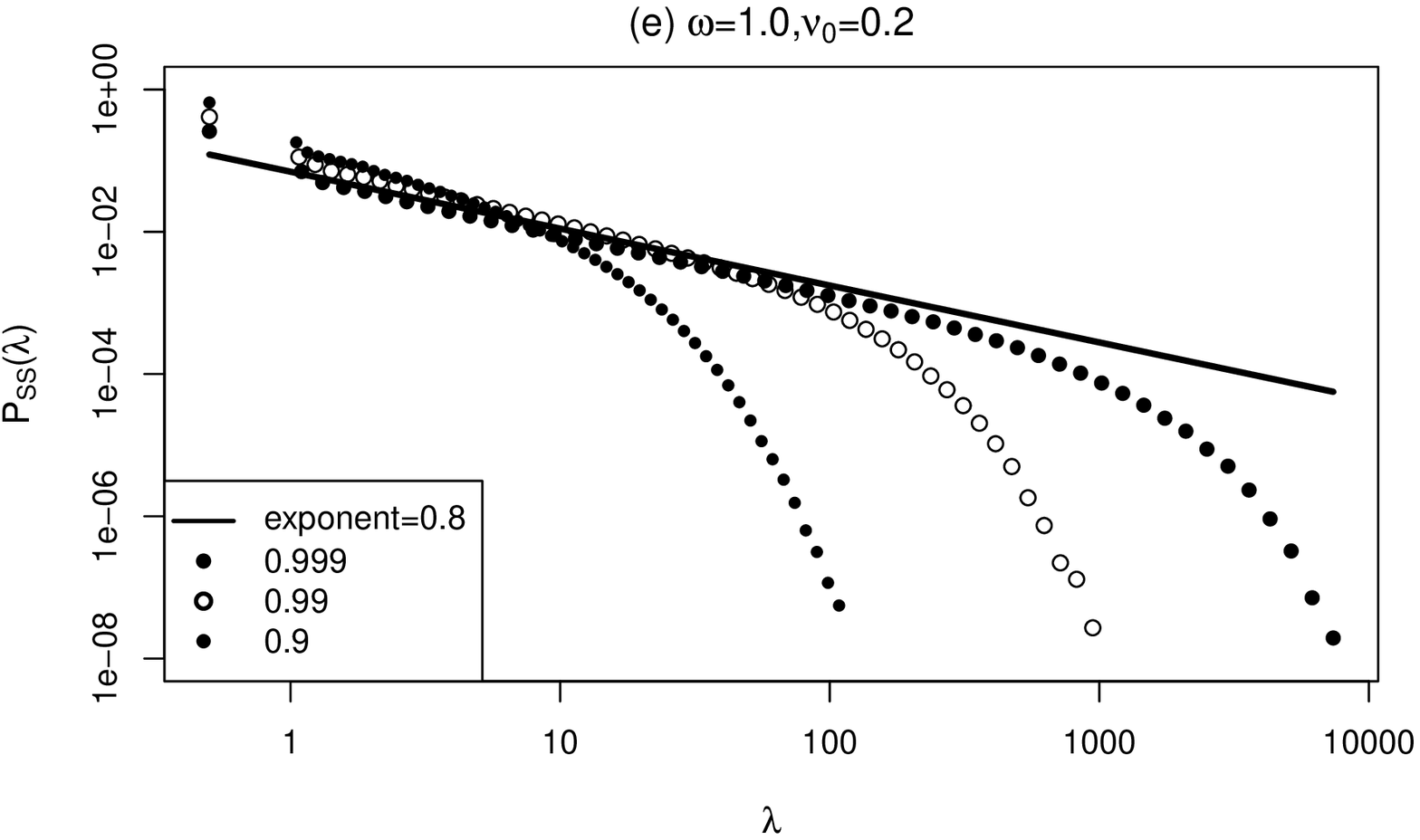}
&   
\includegraphics[width=5cm]{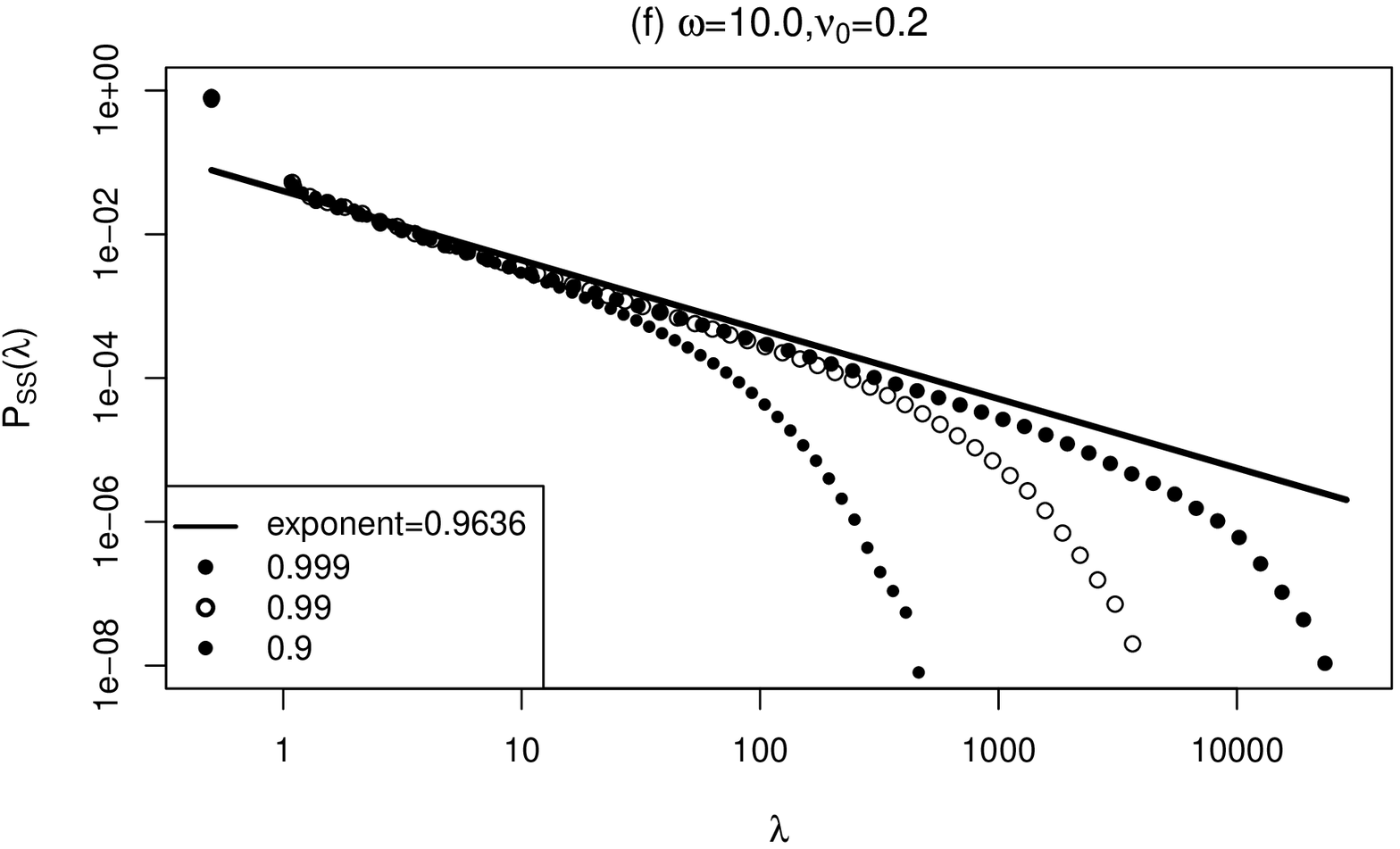}
\\
\includegraphics[width=5cm]{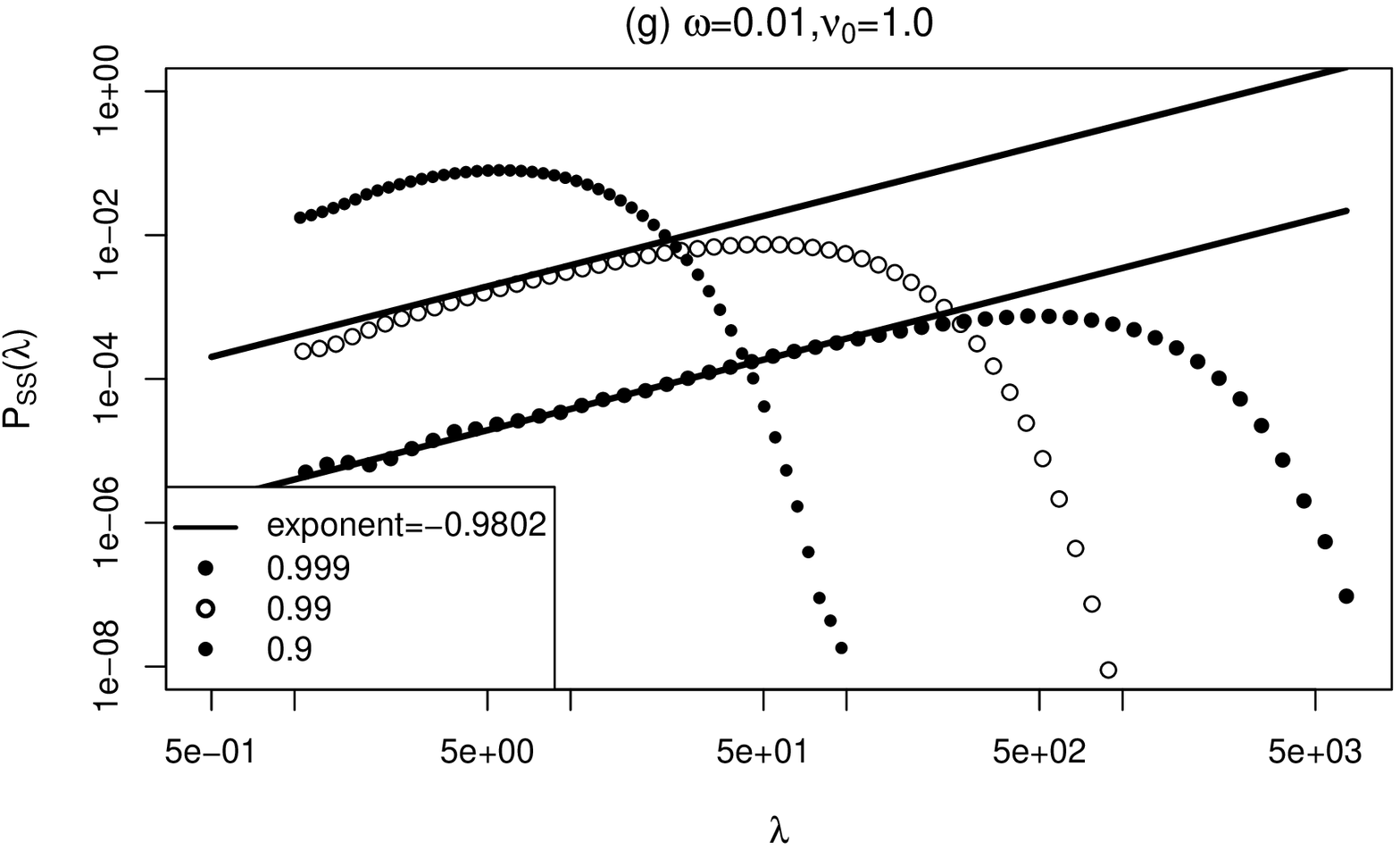}
&   
\includegraphics[width=5cm]{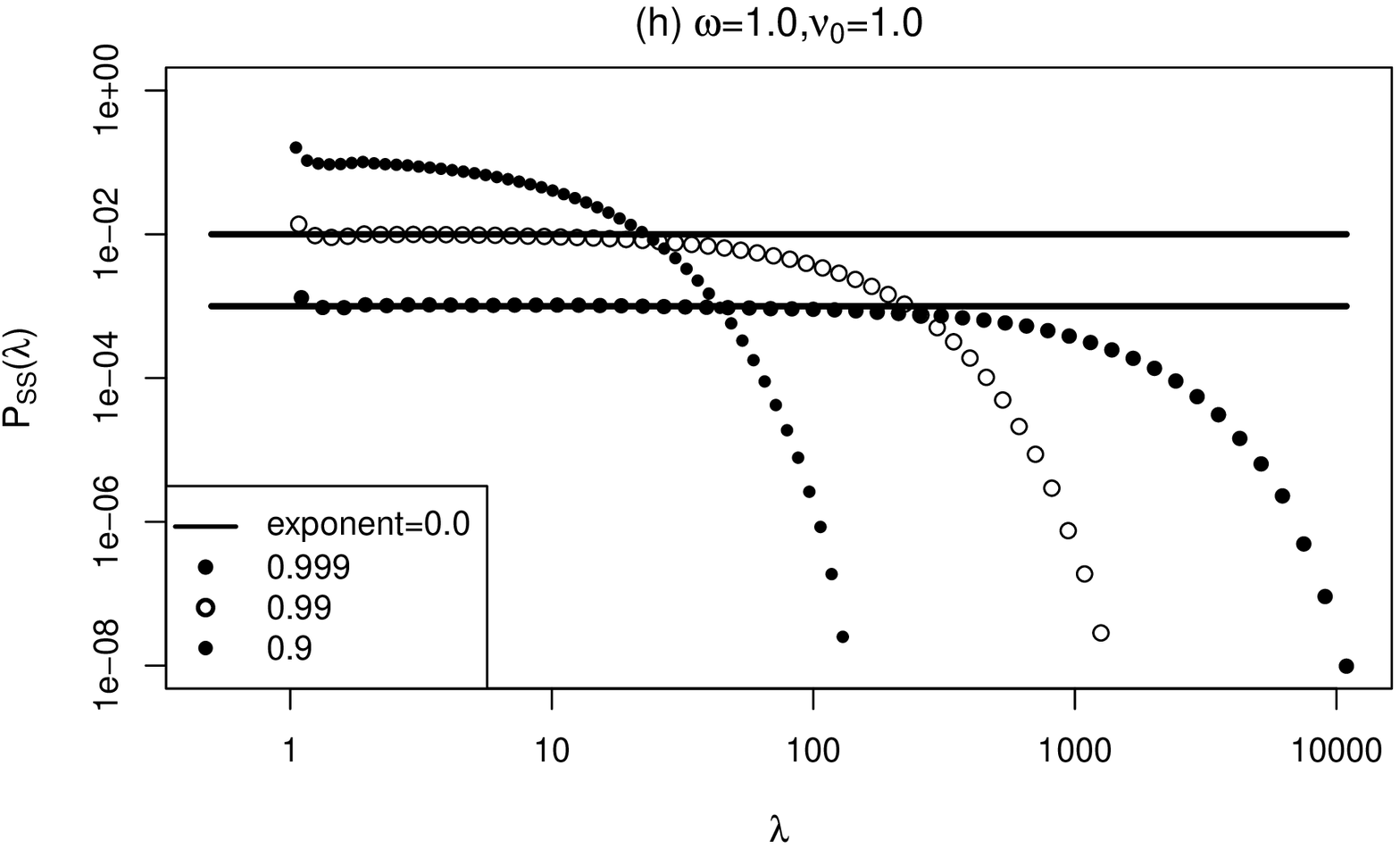}
&   
\includegraphics[width=5cm]{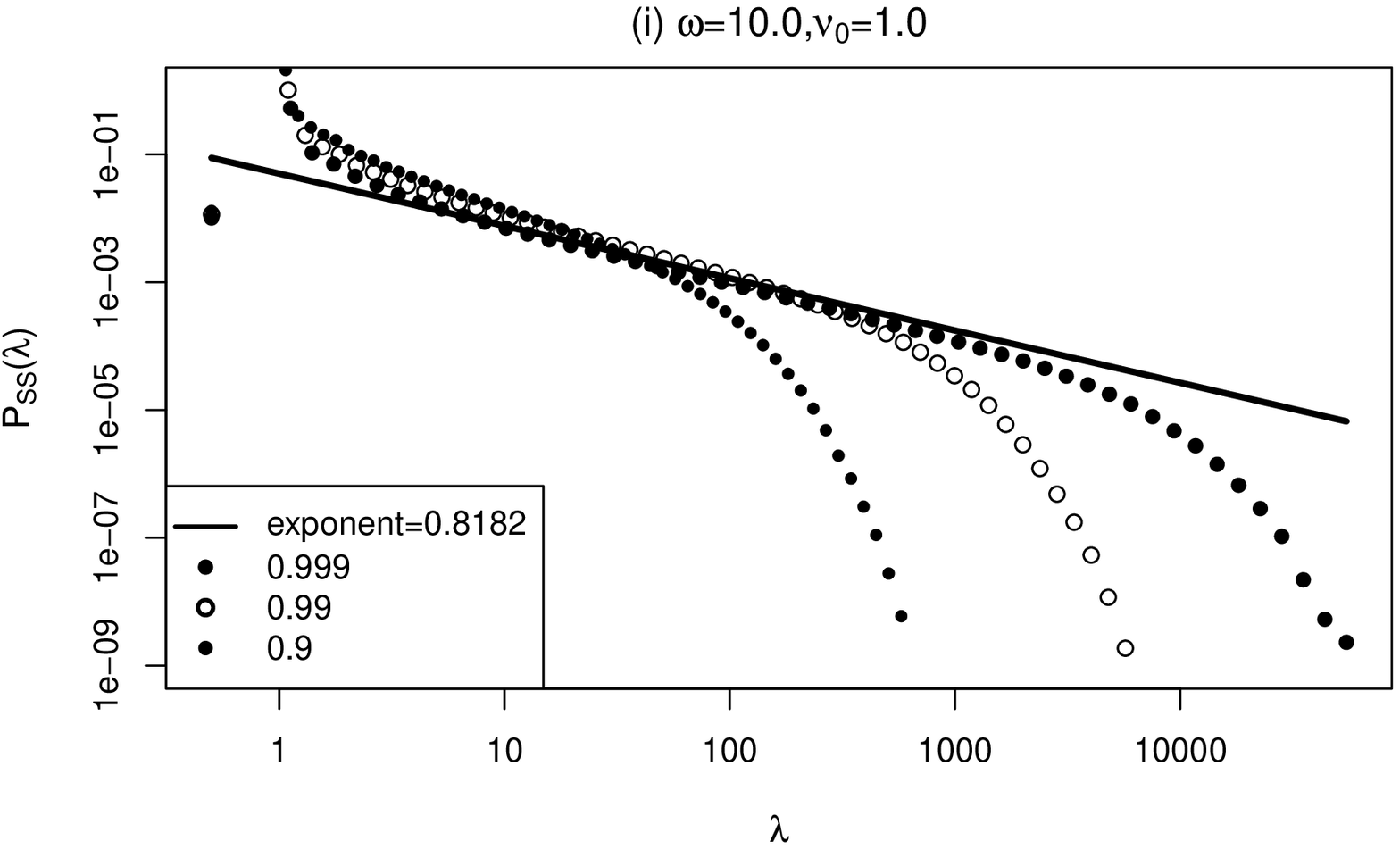}
\end{tabular}
\end{center}
\caption{Plot of the steady state PDF $P_{ss}(\lambda)$ for the single exponential kernel case $K=1$ near the critical point $n=n_c=1$. We set $n=0.9(\bullet)$, $n=0.99)(\circ)$, and $n=0.999)(\bullet)$ for $\tau=1.0,\omega\in\{0.01,1.0,10.0\}$ and $\nu_0\in\{0.01,0.2,1.0\}$.
  (a-c) Background intensity $\nu_0=0.01$ and power-law exponents are $0.9802$, $0.99$, and $0.9982$ for $\omega=0.01$, $1.0$, and $10$, respectively.
  (d-f) Background intensity $\nu_0=0.2$ and power-law exponents are $0.6039$, $0.8$, and $0.9636$ for $\omega=0.01$, $1.0$, and $10$, respectively.
  (g-i) Background intensity $\nu_0=1.0$ and power-law exponents are $-0.9802$, $0.0$, and $0.8182$ for $\omega=0.01$, $1.0$, and $10$, respectively.
} 
\label{fig:dist_lambda_1}
\end{figure}

Here, we test the theoretical prediction of Eq. \eqref{power-law_exponent:1},
\[
P_{SS}(\lambda)\propto \lambda^{-1+2\frac{\nu_0 \tau}{\omega +1}}e^{-\frac{2\tau\epsilon}{\omega +1}\lambda}.
\]
Fig.\ref{fig:dist_lambda_1} presents the PDFs of $\hat{\lambda}(t)$ for the cases.

For a small background intensity $\nu_0=0.01 \ll 1$, the power-law
exponents are $0.9802$ and $0.9982$ for $\omega=0.01$ and $\omega=10.0$,
respectively. These exponents are close to 1 and show little dependence on $\omega$. 
For large $\nu_0$, the power-law exponent becomes more dependent on $\omega$.
When $\nu_0=1.0$, the power-law exponents are $-0.9802$ and $0.8182$ for
$\omega=0.01$ and $\omega=10.0$, respectively.
When $\omega \gg 2\nu_0\tau$,
The power-law exponent becomes universally 1, as the expression
$1-2\nu_0\tau/(\omega+1)\simeq 1$ holds.
Furthermore, the length scale of the power-law region, denoted as $(\omega + 1)/(2\tau\epsilon)$, increases as $\omega$ increases.
The numerical results confirm these findings by demonstrating agreement between the slopes of the PDFs  and the theoretical values, and the exponents of the power-law were verified. In addition, by observing the $x$-axis, it can be seen that the range of the straight power-law region becomes wider as $\omega$ increases.

\subsubsection*{b. Double exponential kernel}
\begin{figure}[htbp]
\begin{center}
\begin{tabular}{ccc}    
\includegraphics[width=5cm]{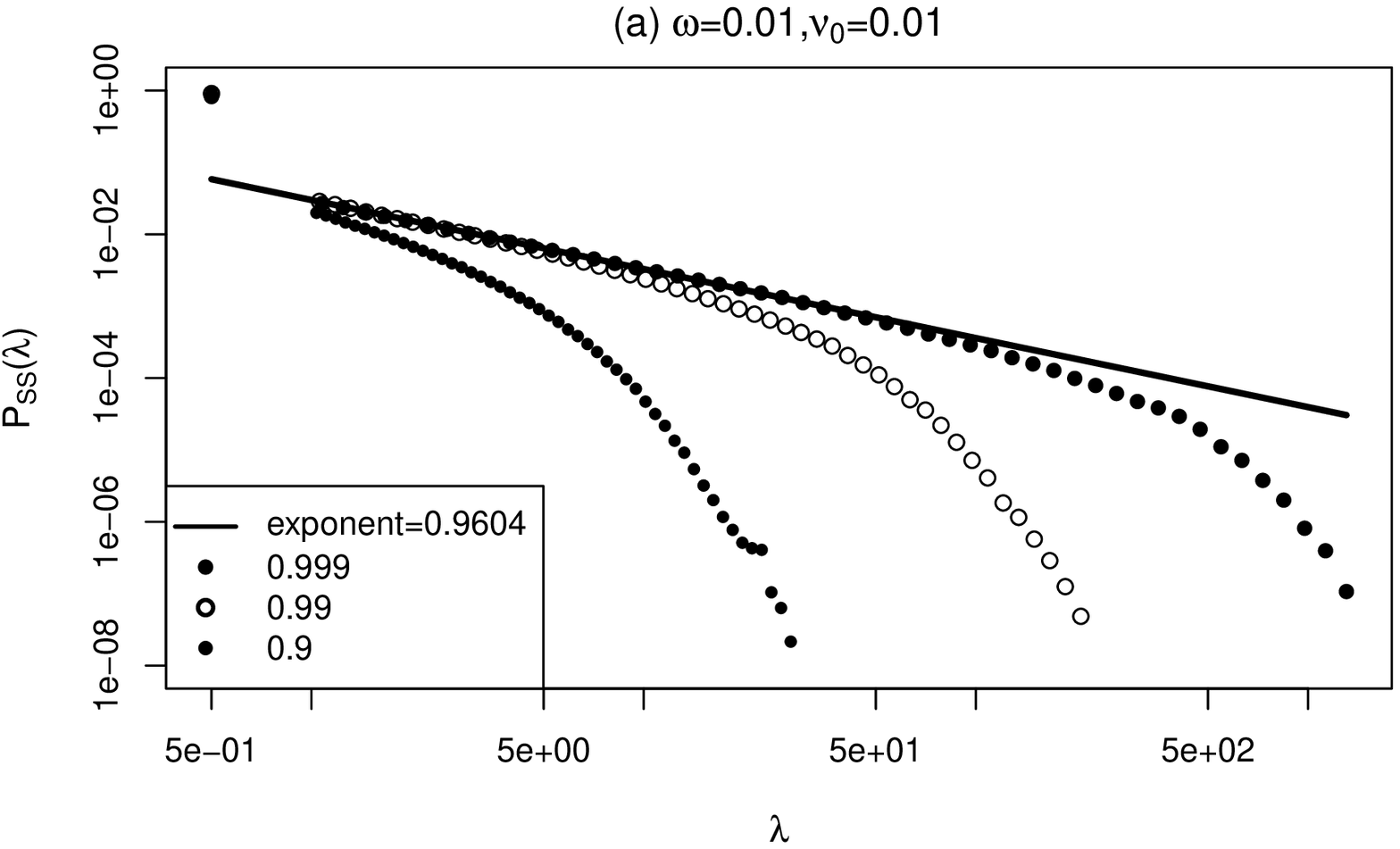}
&   
\includegraphics[width=5cm]{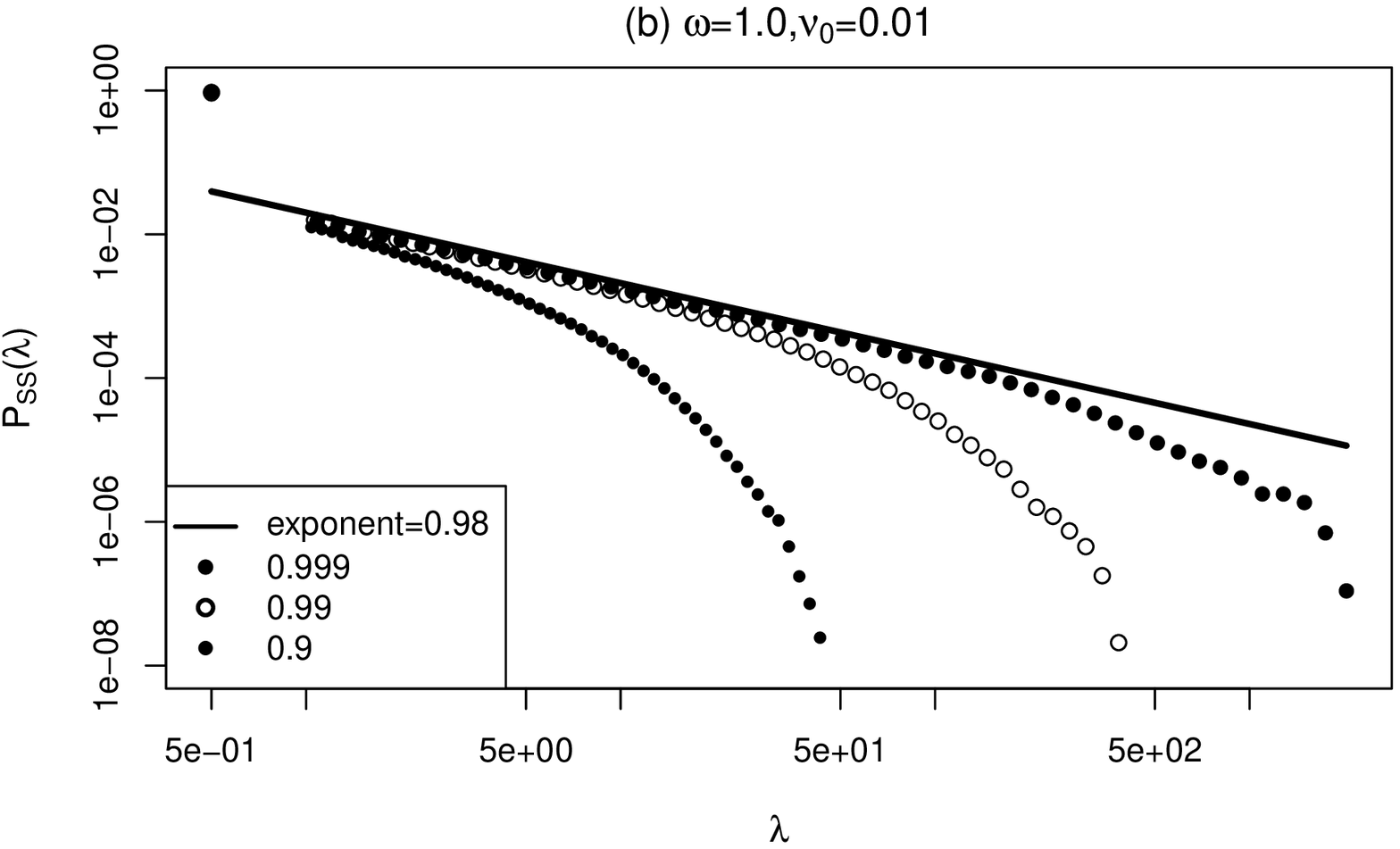}
&   
\includegraphics[width=5cm]{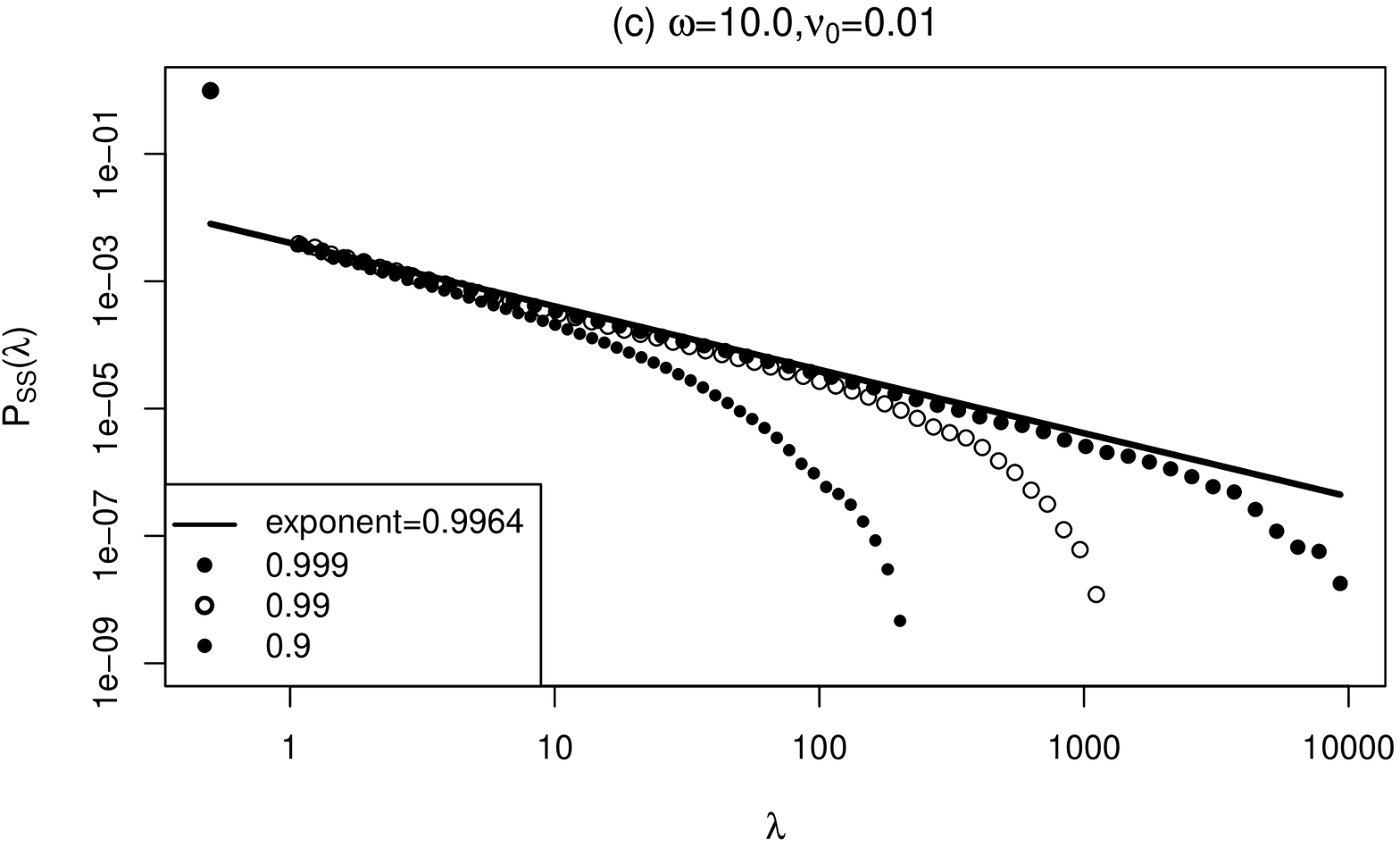}
\\
\includegraphics[width=5cm]{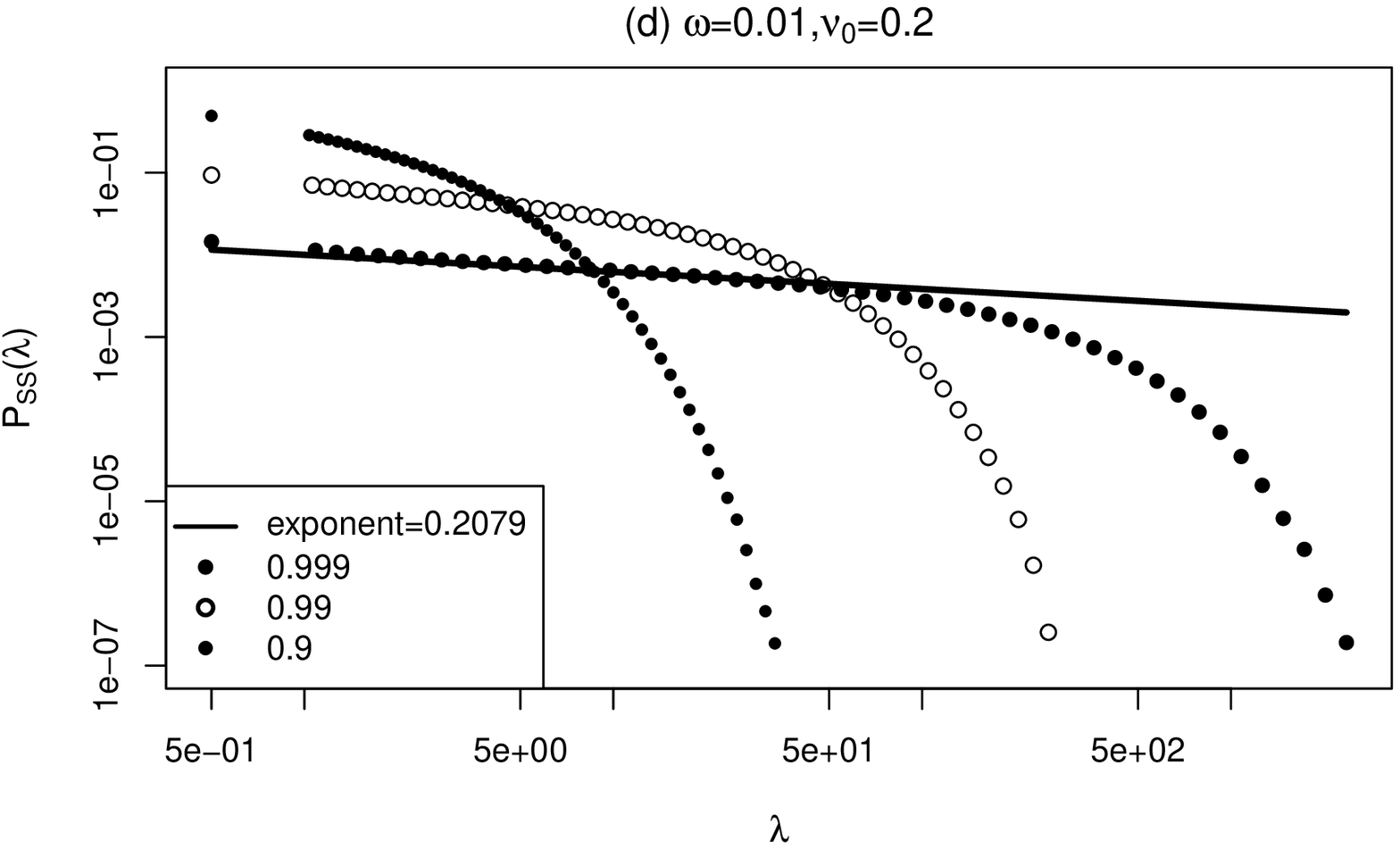}
&   
\includegraphics[width=5cm]{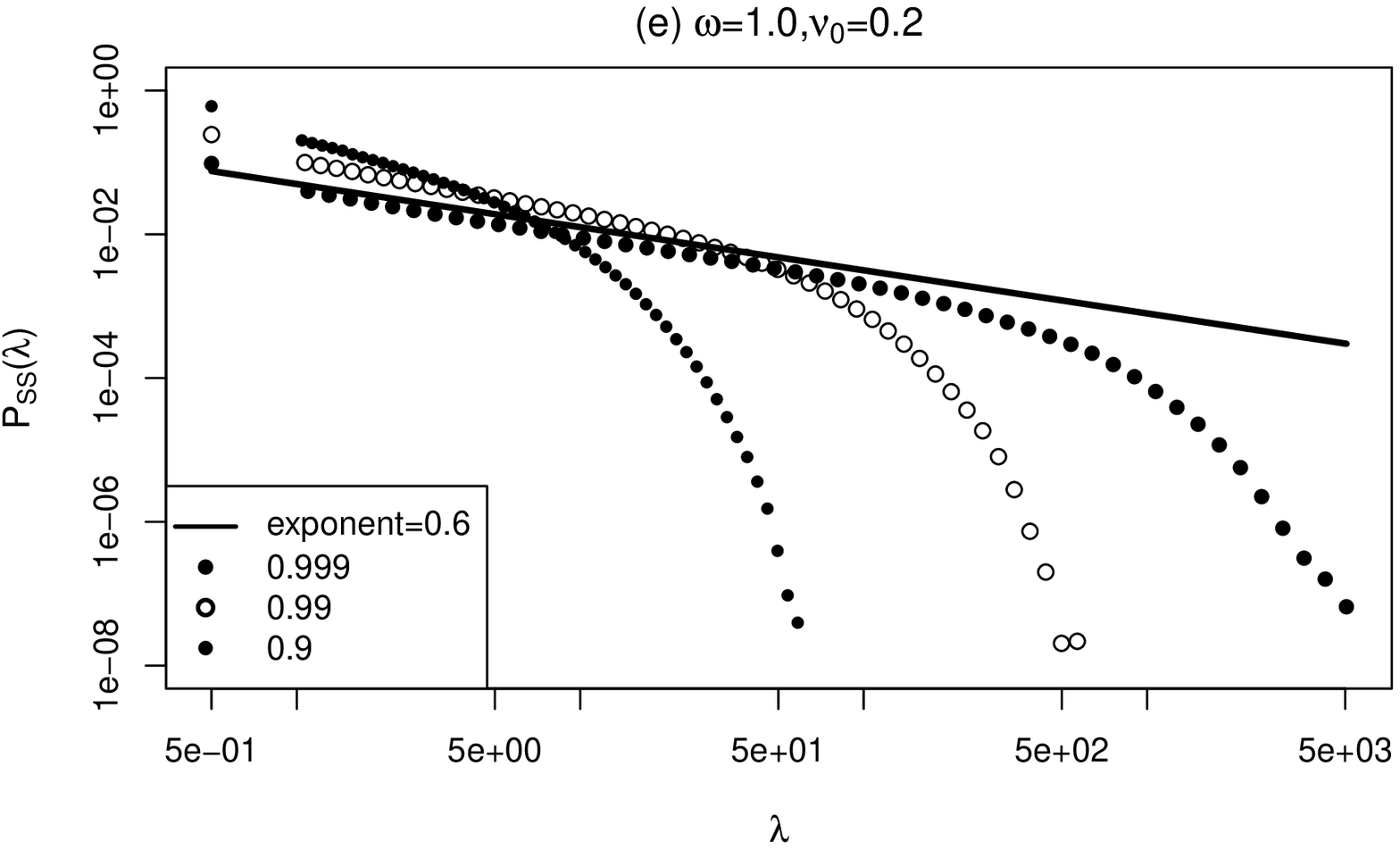}
&   
\includegraphics[width=5cm]{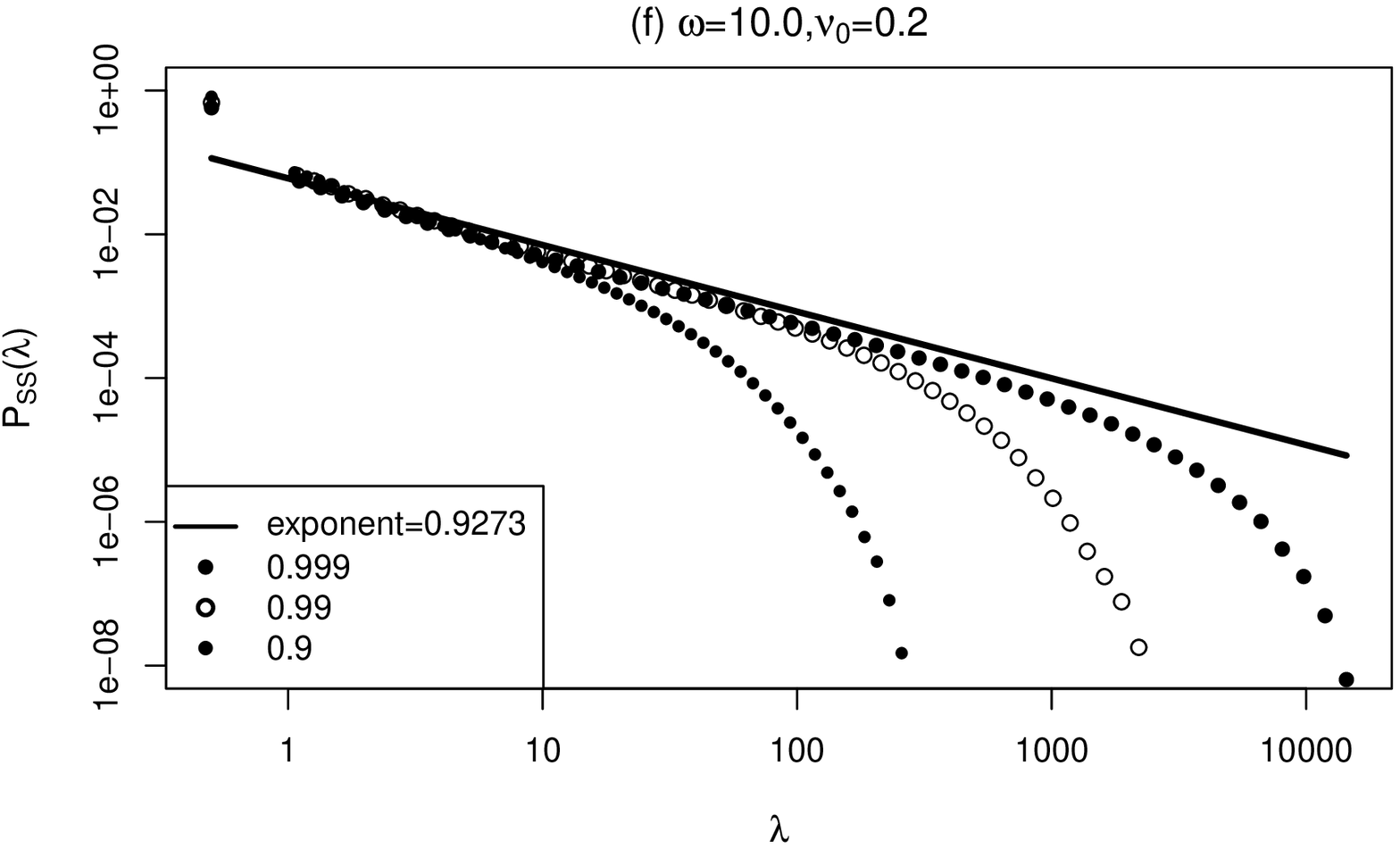}
\\
\includegraphics[width=5cm]{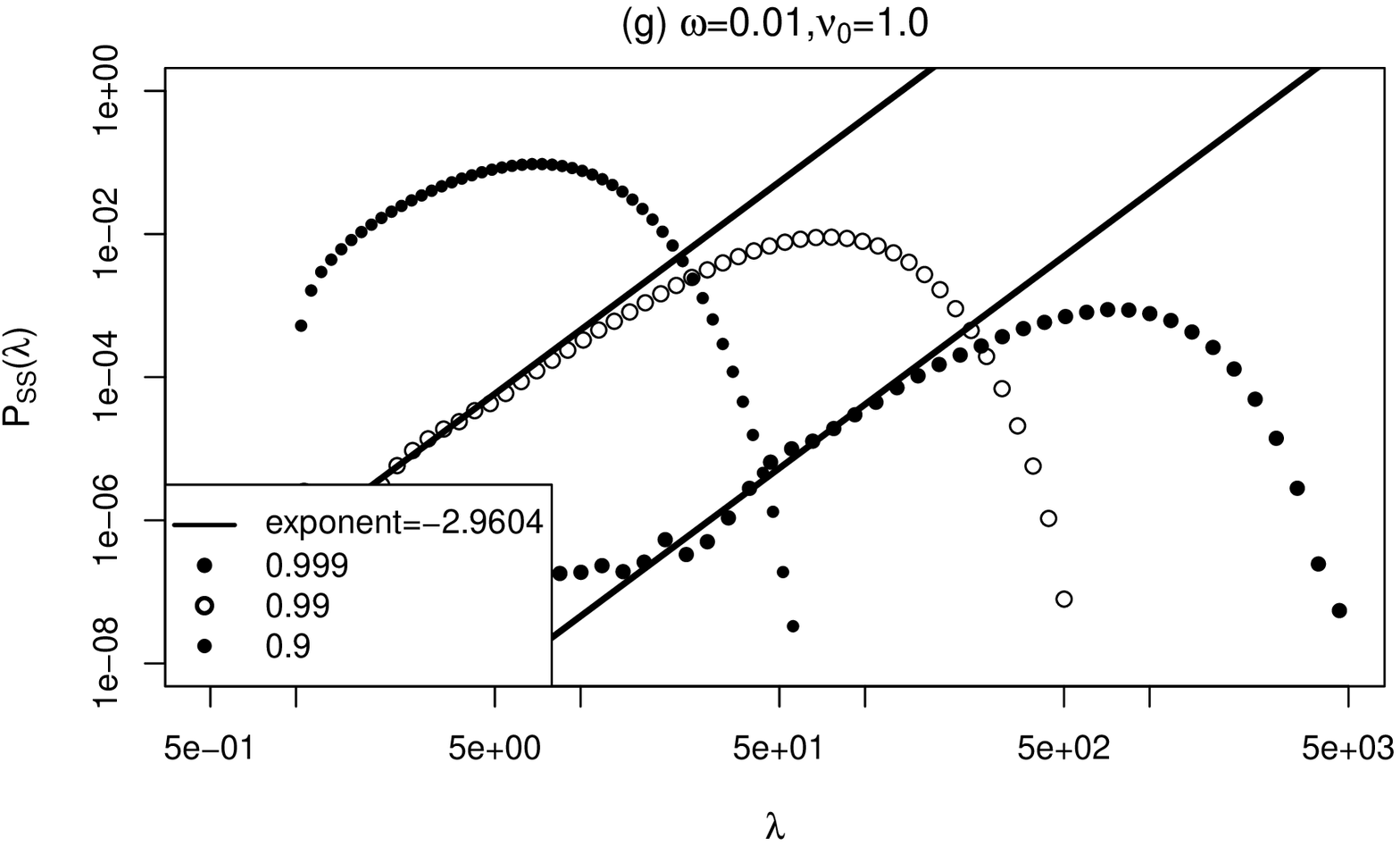}
&   
\includegraphics[width=5cm]{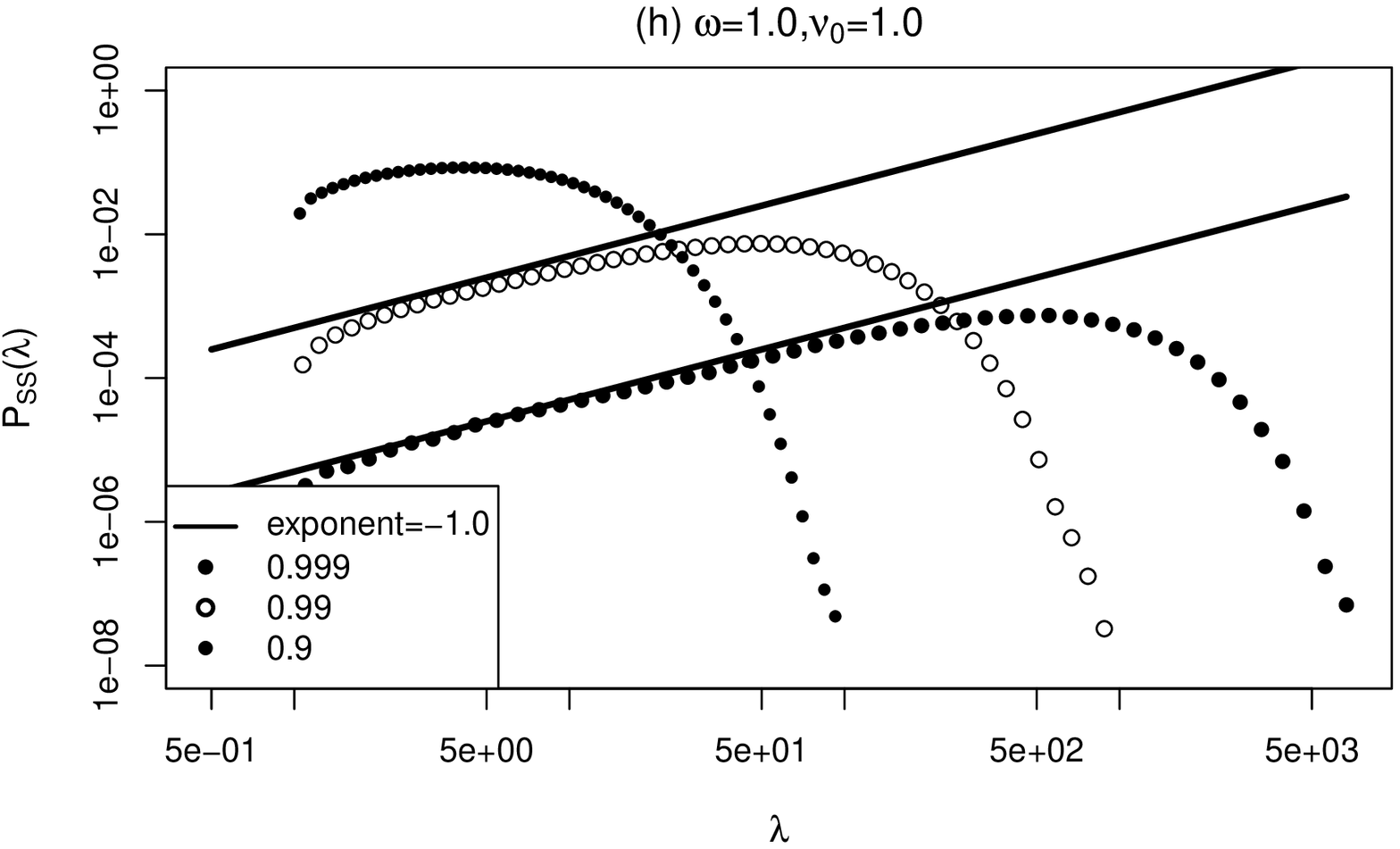}
&   
\includegraphics[width=5cm]{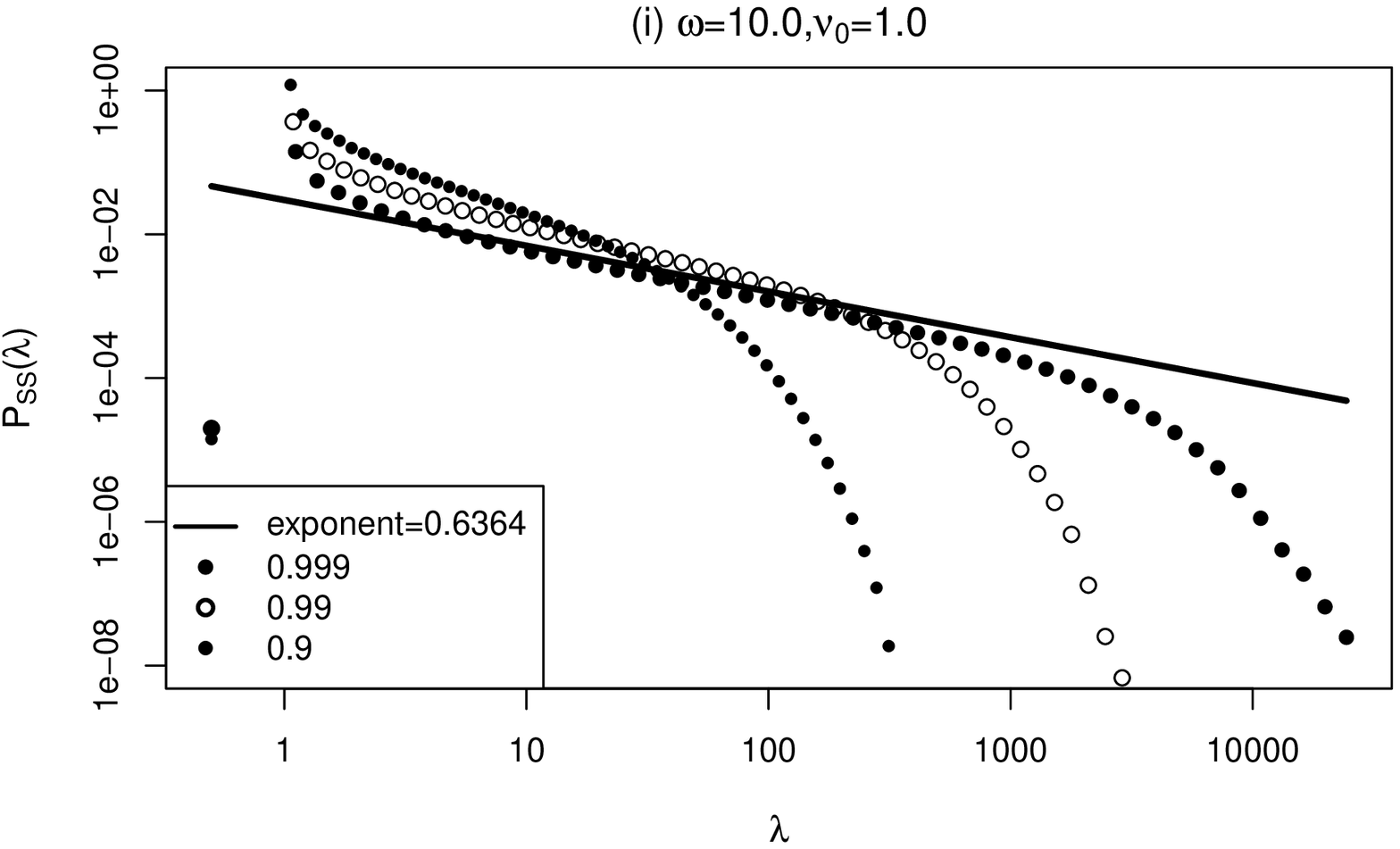}
\end{tabular}
\end{center}
\caption{
Plot of the steady-state PDF $P_{ss}(\lambda)$ for the double exponential kernel case $K=2$ near the critical point $n=n_c=1$. We set $n=0.0(\bullet)$, $n=0.99(\circ)$, and $n=0.999(\bullet)$ for $(\tau_1,\tau_2)=(1.0,3.0),\omega\in\{0.01,1.0,10.0\}$, and $\nu_0\in\{0.01,0.2,1.0\}$.
  (a-c) Background intensity $\nu_0=0.01$ and power-law exponents are $0.9604$, $0.98$, and $0.9964$ for $\omega=0.01$, $1.0$, and $10$, respectively.
  (d-f) Background intensity $\nu_0=0.2$ and power-law exponents are $0.2079$, $0.6$, and $0.9273$ for $\omega=0.01$, $1.0$, and $10$, respectively.
  (g-i) Background intensity $\nu_0=1.0$ and power-law exponents are $-2.96$, $-1.0$, and $0.6363$ for $\omega=0.01$, $1.0$, and $10$, respectively.
} 
\label{fig:dist_lambda_2}
\end{figure}

We verified the theoretical predictions of Eq. \eqref{power-law_exponent:2}.
\[
P_{SS}(\lambda)\propto \lambda^{-1+2\frac{\nu_0 (\tau_1 n_1+\tau_2 n_2)}{\omega +1}}.
\]
To realize $n=n_1+n_2\in \{0.9,0.99,0.999\}$, we set
$(n_1,n_2)=(0.5,0.4), (0.5,0.49)$ and $(0.5,0.499)$.
Fig.\ref{fig:dist_lambda_2} presents the double logarithmic plot of 
PDFs of $\hat{\lambda}_t$.
One can observe the agreement between the slopes of the PDFs  and the theoretical values 
and the theoretical prediction is verified.

\subsubsection*{c. Triple exponential kernel}
\begin{figure}[htbp]
\begin{center}
\begin{tabular}{ccc}
\includegraphics[width=5cm]{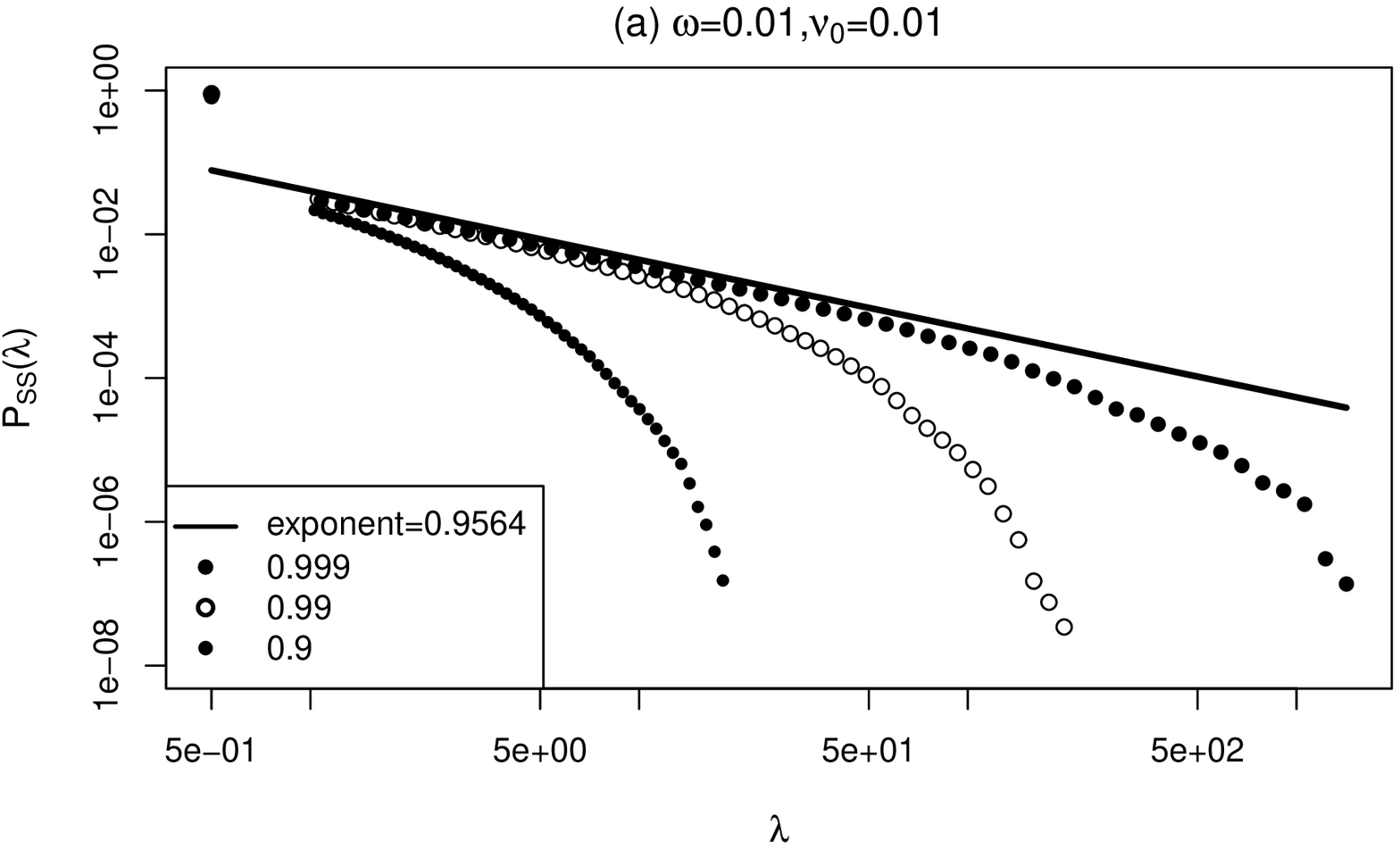}
&   
\includegraphics[width=5cm]{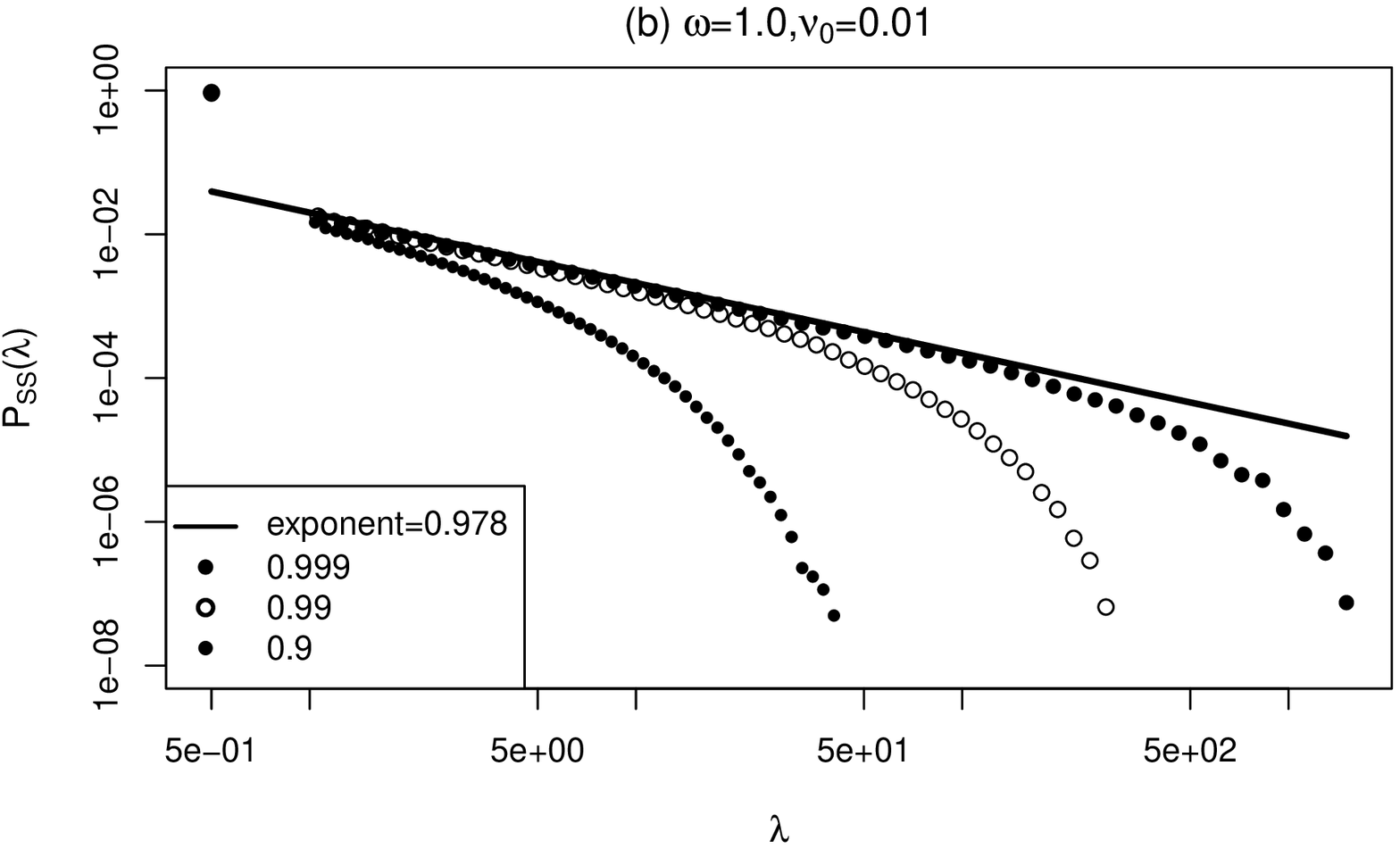}
&   
\includegraphics[width=5cm]{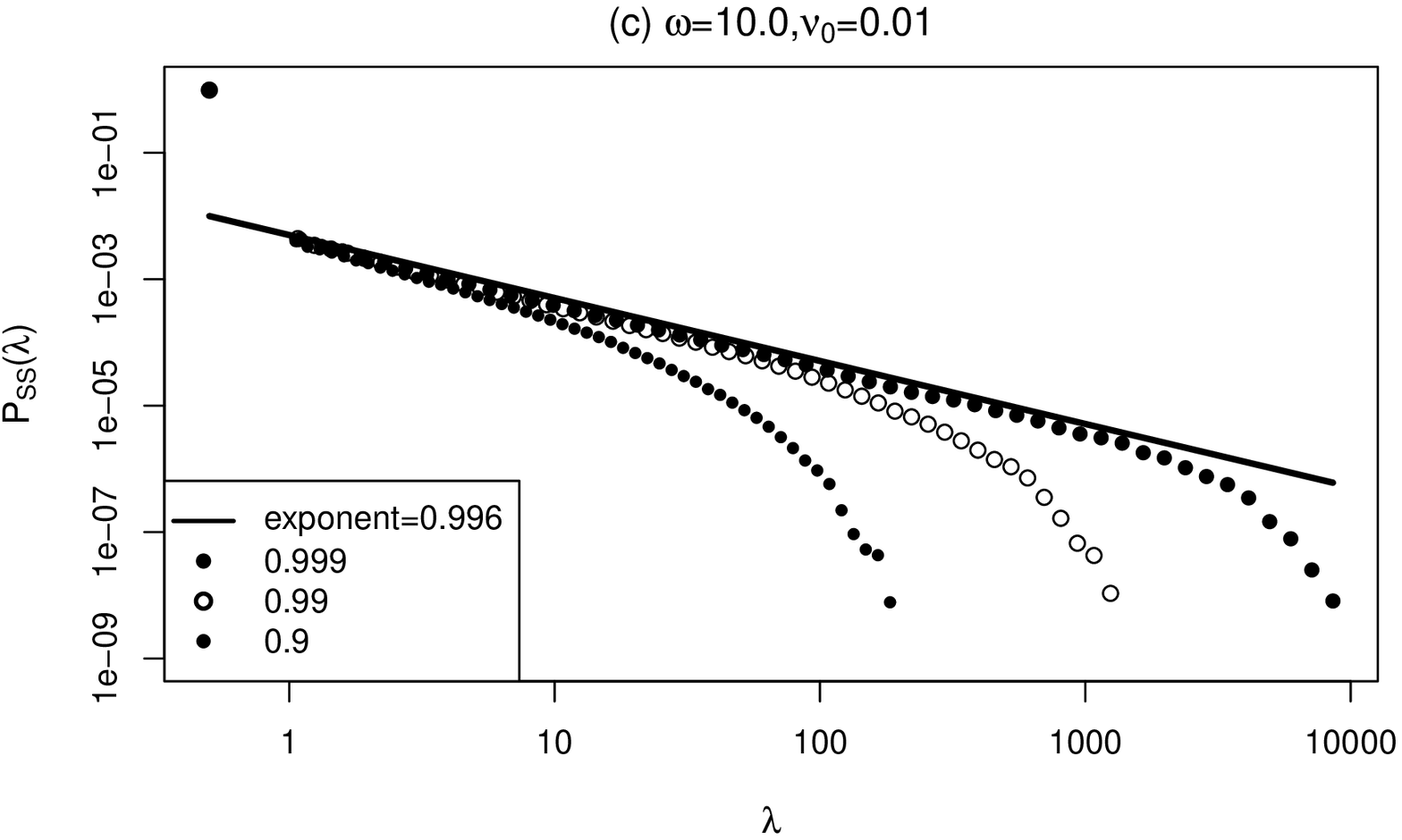}
\\
\includegraphics[width=5cm]{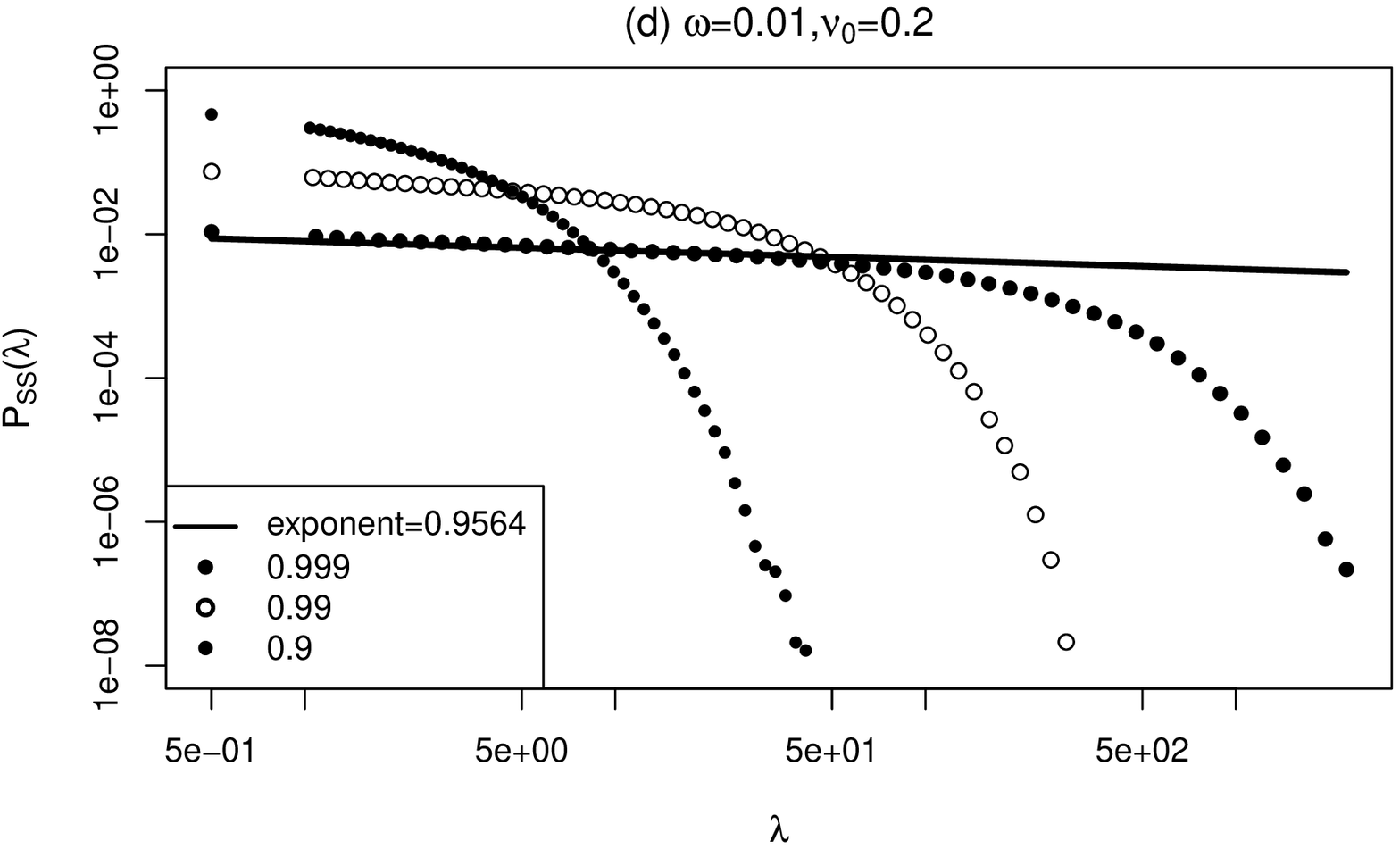}
&   
\includegraphics[width=5cm]{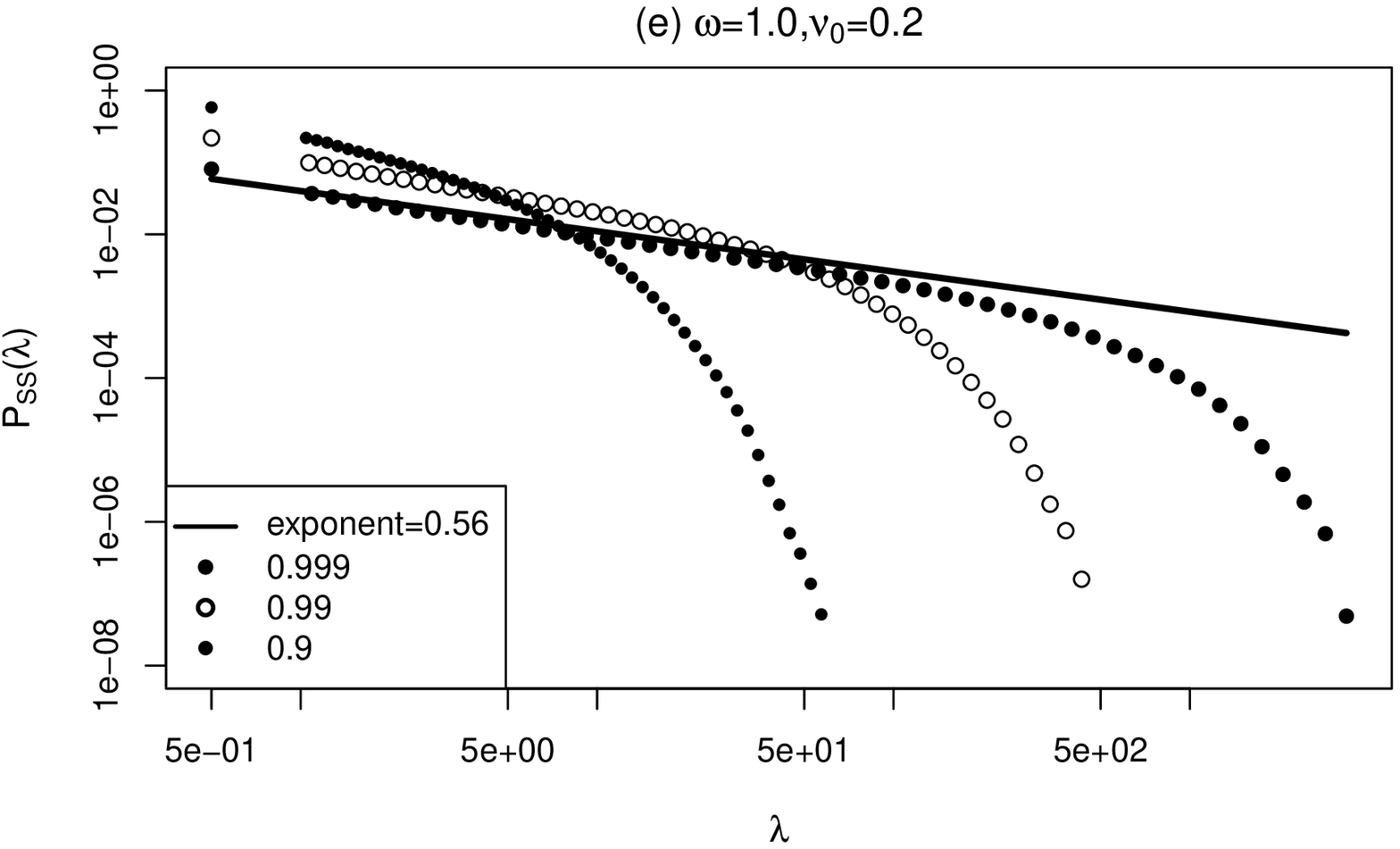}
&   
\includegraphics[width=5cm]{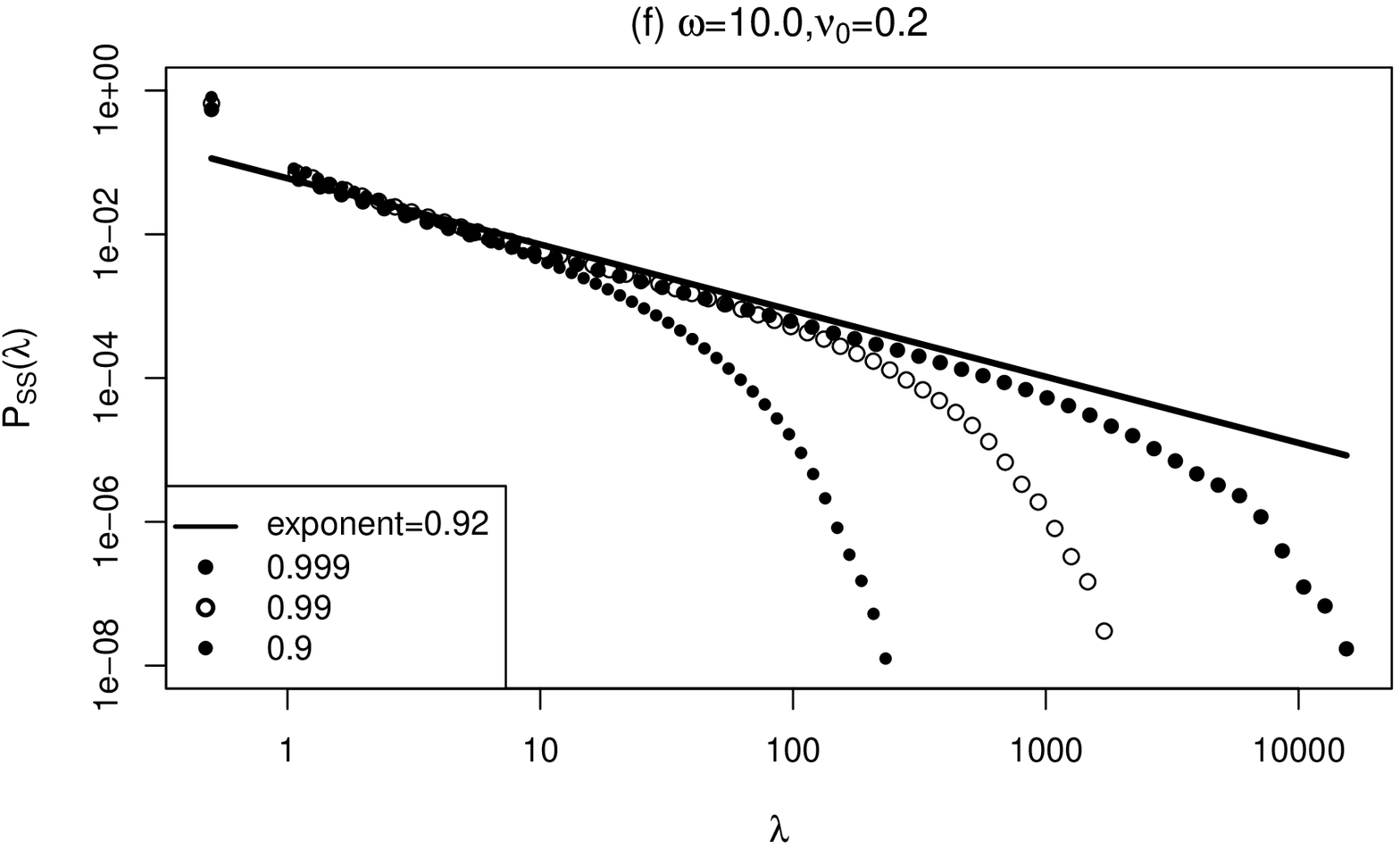}
\\
\includegraphics[width=5cm]{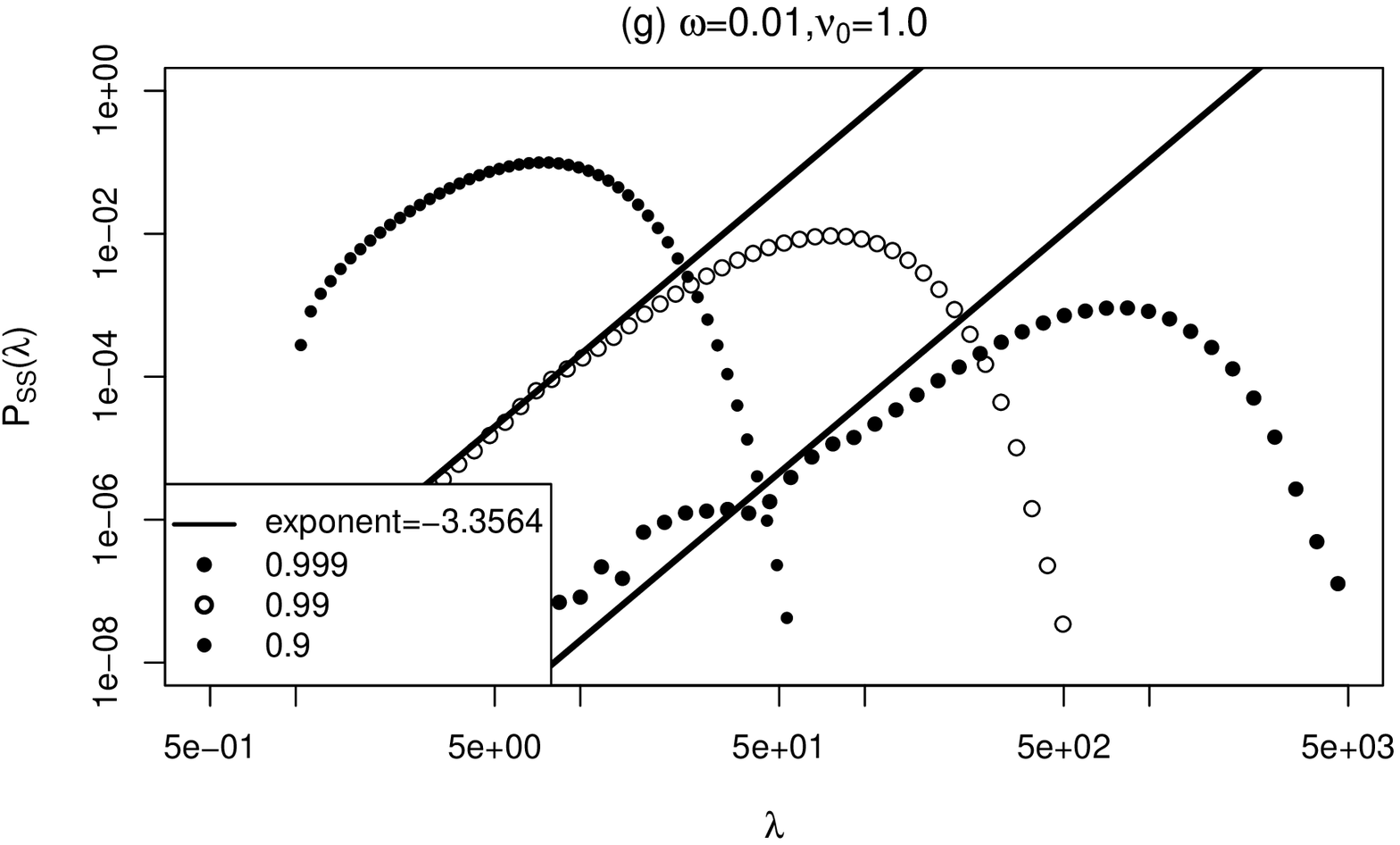}
&   
\includegraphics[width=5cm]{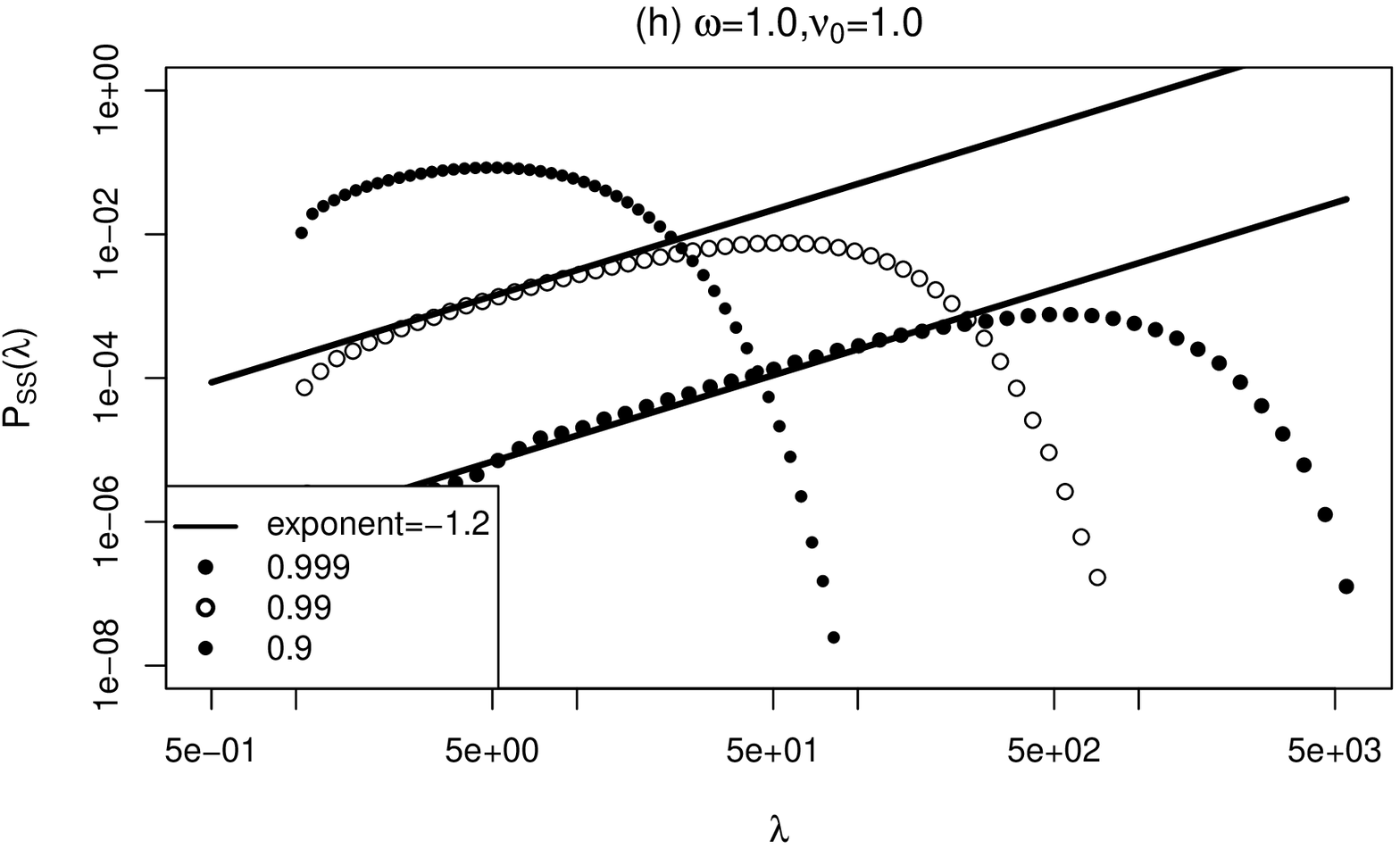}
&   
\includegraphics[width=5cm]{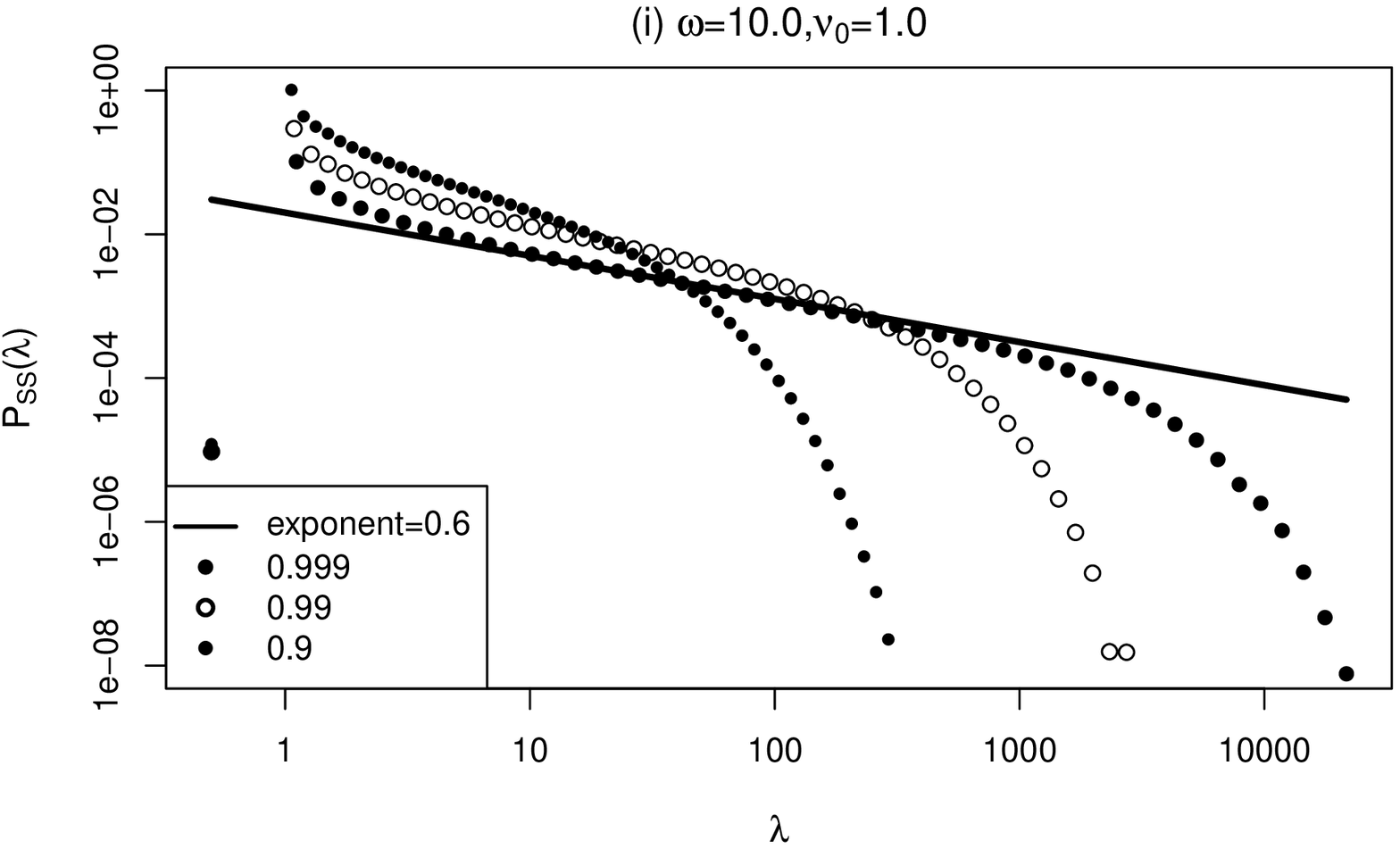}
\end{tabular}
\end{center}
\caption{
  Plot of the numerical results for the steady-state PDF $P_{ss}(\lambda)$ for the double exponential case $K=3$. $\lambda$ near the critical point $n=n_c=1$. We set $n=0.9(\bullet)$, $n=0.99(\circ)$, and $n=0.999(\bullet)$ for $(\tau_1,\tau_2,\tau_3)=(1.0,2.0,3.0),\omega\in\{0.01,1.0,10.0\}$, and $\nu_0\in\{0.01,0.2,1.0\}$.
  (a-c) Background intensity $\nu_0=0.01$ and power-law exponents are $0.956$, $0.978$, and $0.996$ for $\omega=0.01$, $1.0$, and $10$, respectively.
  (d-f) Background intensity $\nu_0=0.2$ and power-law exponents are $0.129$, $0.56$, and $0.92$ for $\omega=0.01$, $1.0$, and $10$, respectively.
  (g-i) Background intensity $\nu_0=1.0$ and power-law exponents are $-3.356$, $-1.2$, and $0.6$ for $\omega=0.01$, $1.0$, and $10$, respectively.
}
\label{fig:dist_lambda_3}
\end{figure}

We verified the theoretical predictions of Eq. \eqref{power-law_exponent:K}
 for $K=3$.
\[
P_{SS}(\lambda)\propto \lambda^{-1+2\frac{\nu_0 (\tau_1 n_1+\tau_2 n_2+\tau_3 n_3)}{\omega +1}}.
\]
To achieve $n=n_1+n_2+n_3\in \{0.9,0.99,0.999\}$, we set
$(n_1,n_2,n_3)=(0.3,0.2,0.4), (0.3,0.2,0.49)$, and $(0.3,0.2,0.499)$.
Fig.\ref{fig:dist_lambda_3} shows the resulting PDFs of $\hat{\lambda}(t)$.
One can observe the agreement between the slopes of the PDFs  and the theoretical values 
and the theoretical prediction is verified.

\subsubsection*{d. Power-law memory kernel}

\begin{figure}[htbp]
\begin{center}
\begin{tabular}{ccc}
\includegraphics[width=5cm]{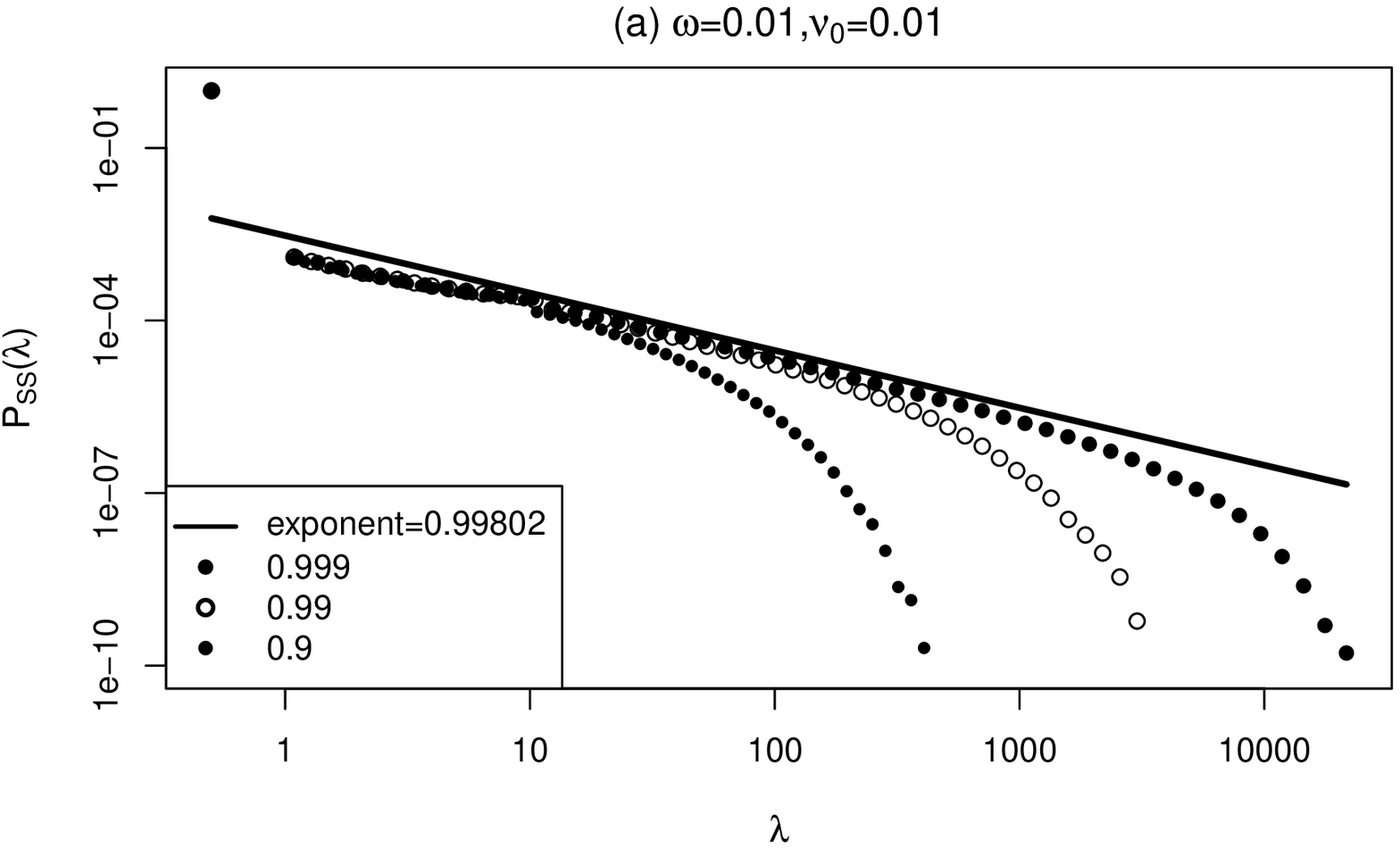}
&   
\includegraphics[width=5cm]{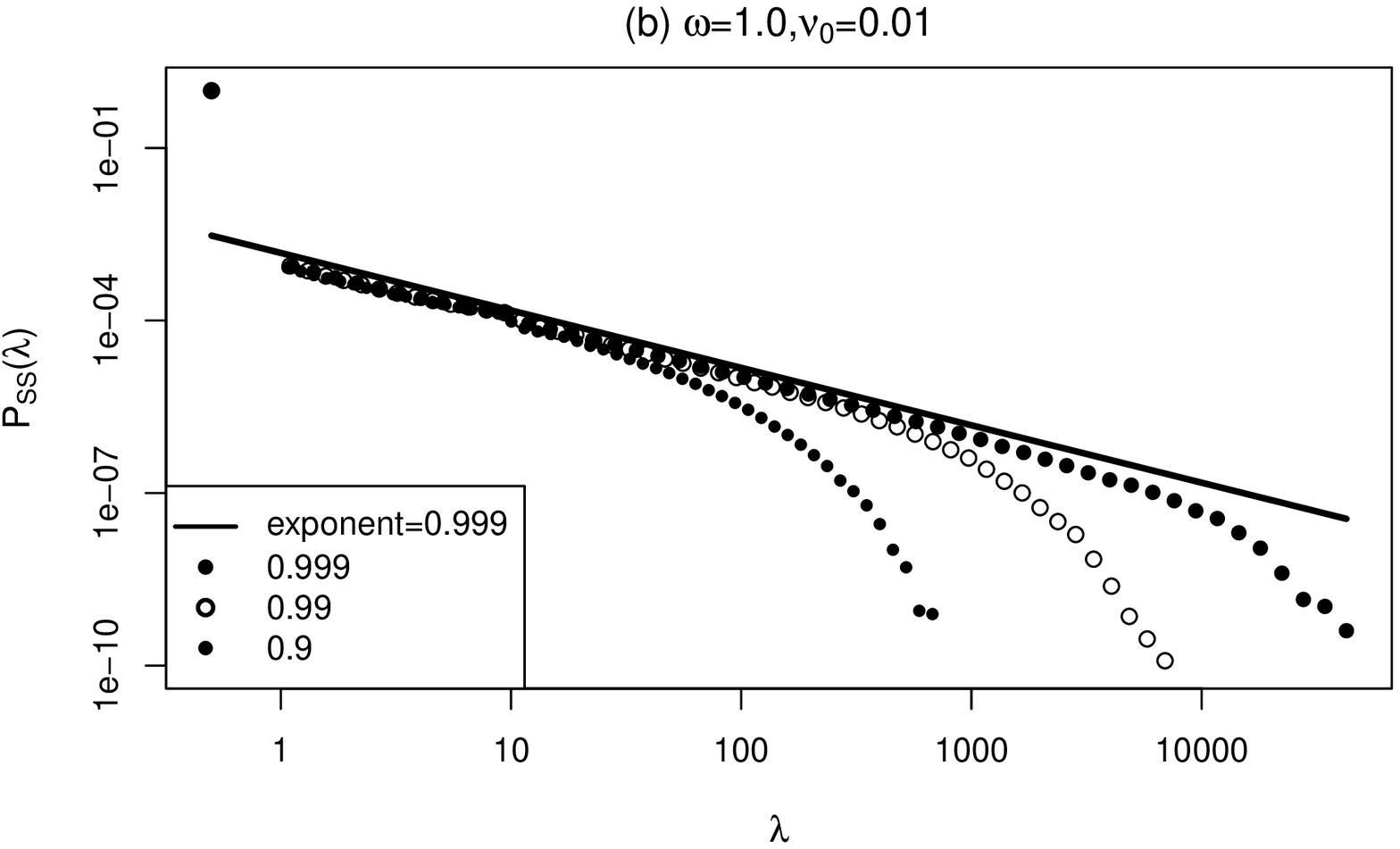}
&   
\includegraphics[width=5cm]{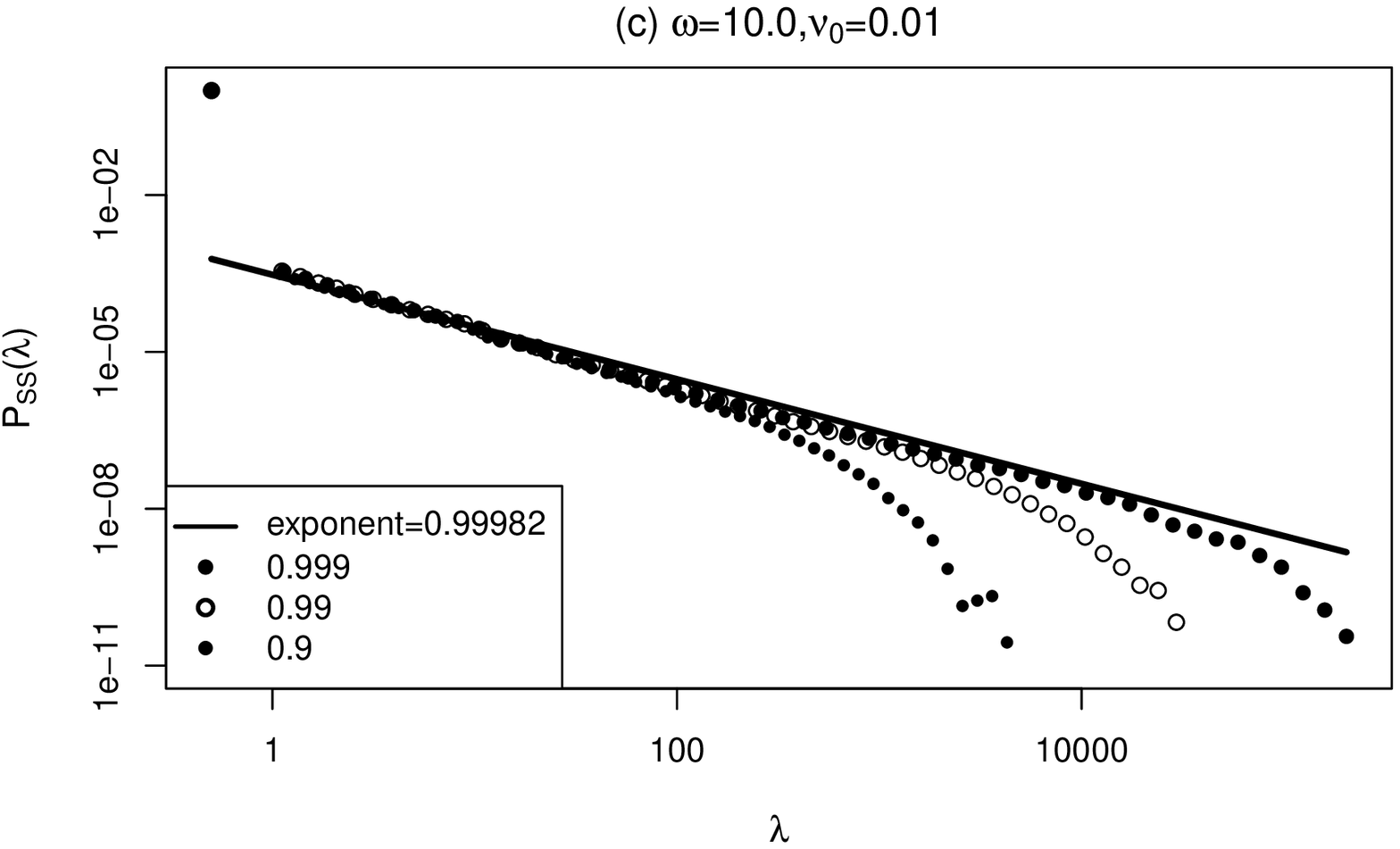}
\\
\includegraphics[width=5cm]{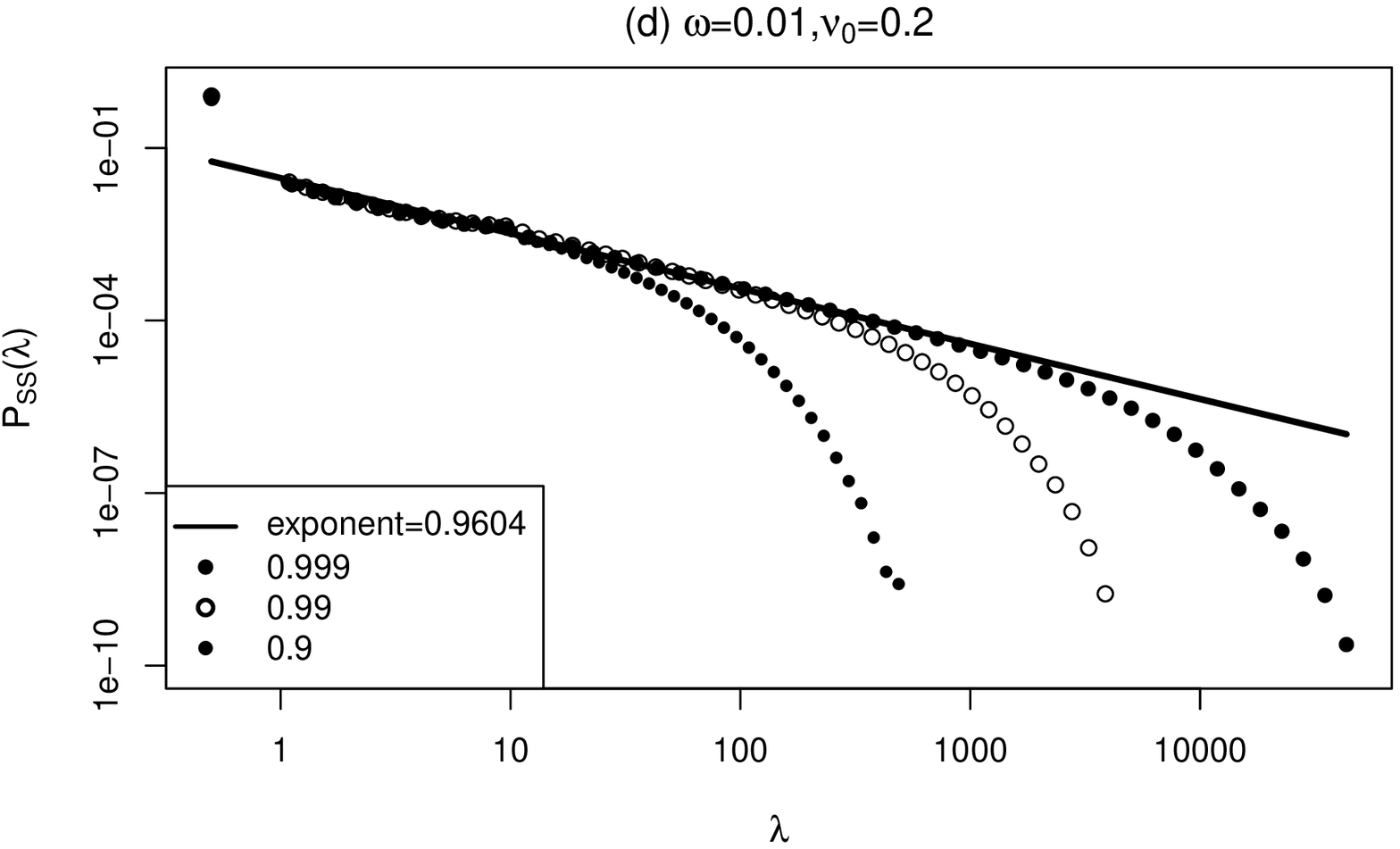}
&   
\includegraphics[width=5cm]{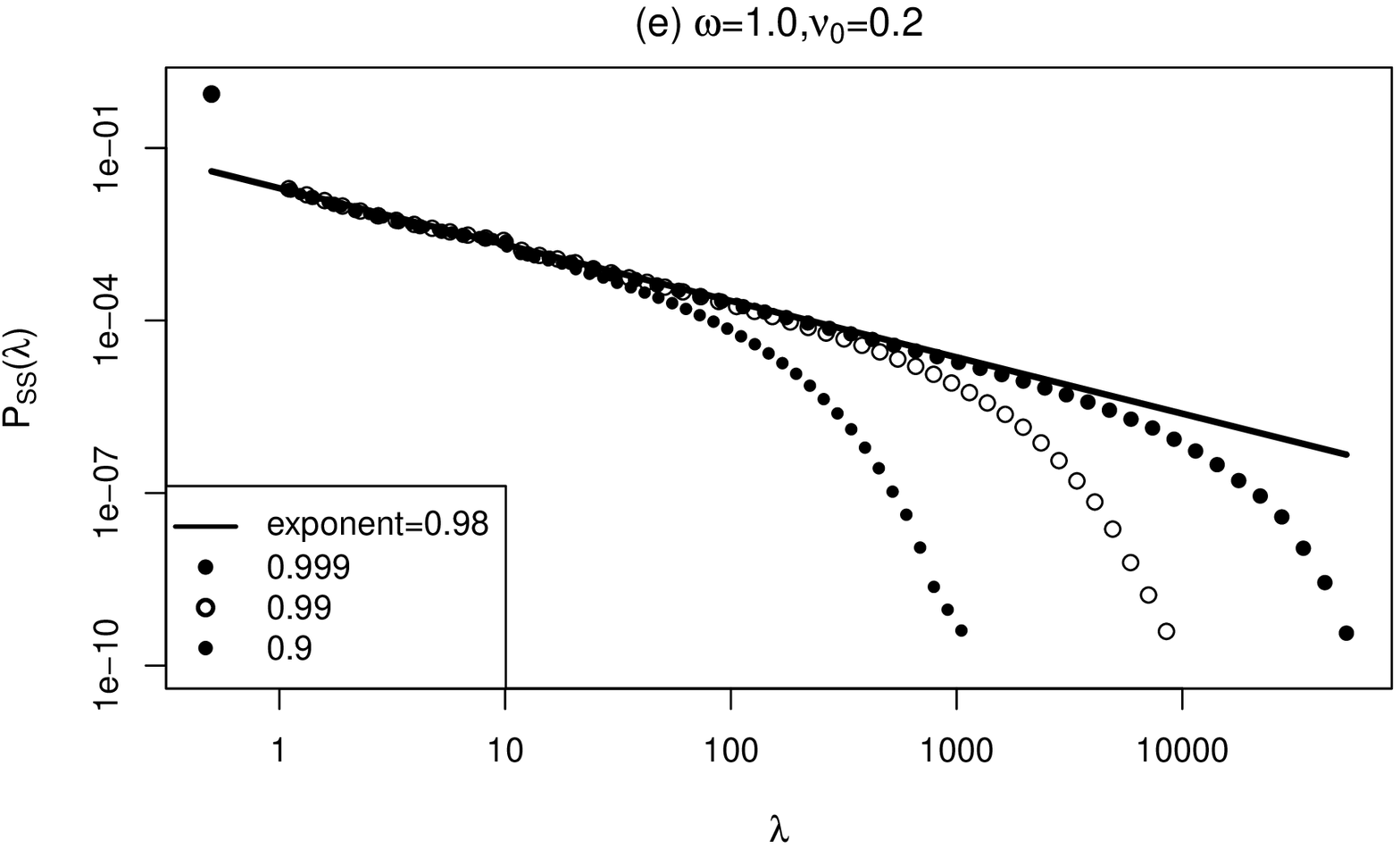}
&   
\includegraphics[width=5cm]{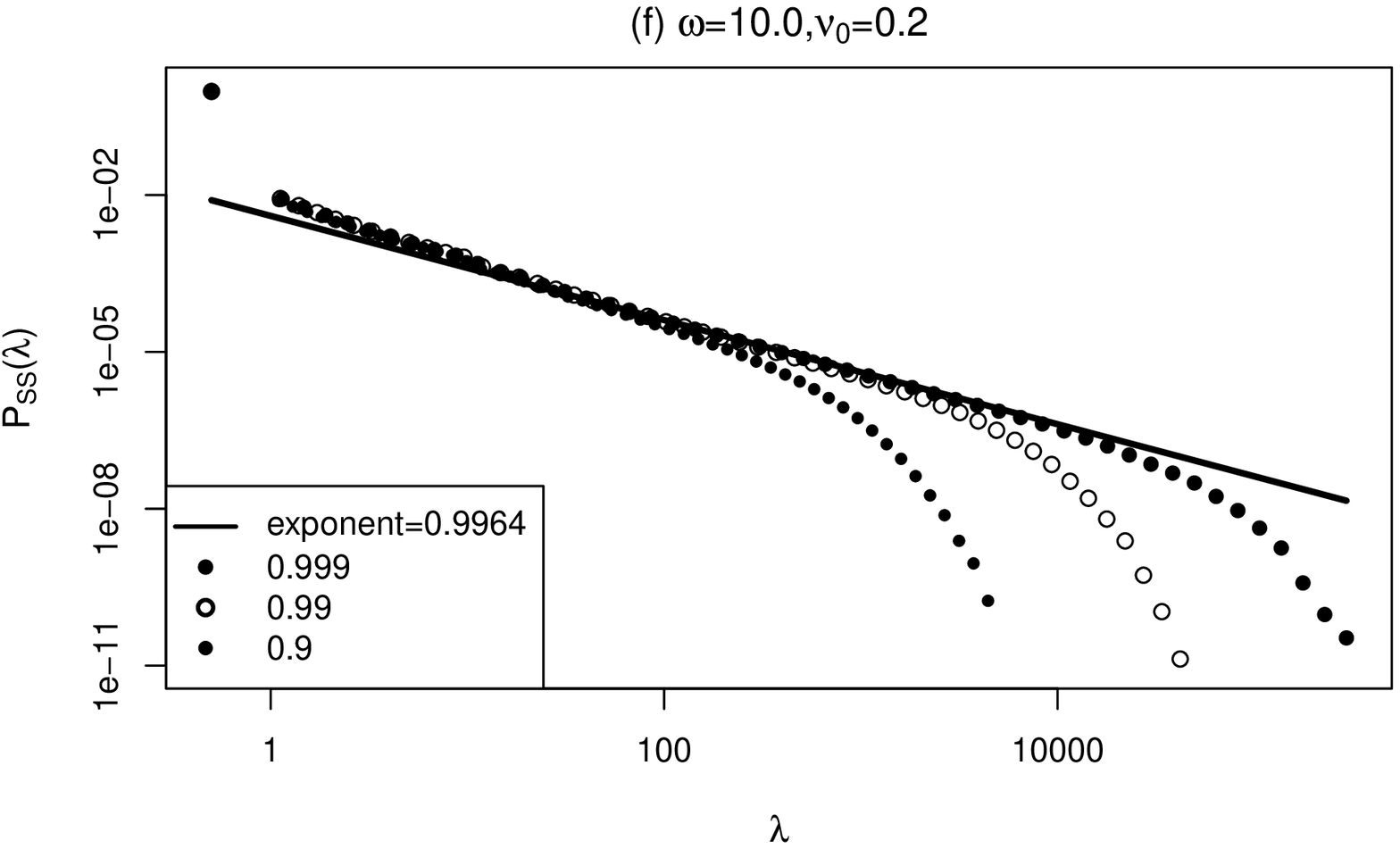}
\\
\includegraphics[width=5cm]{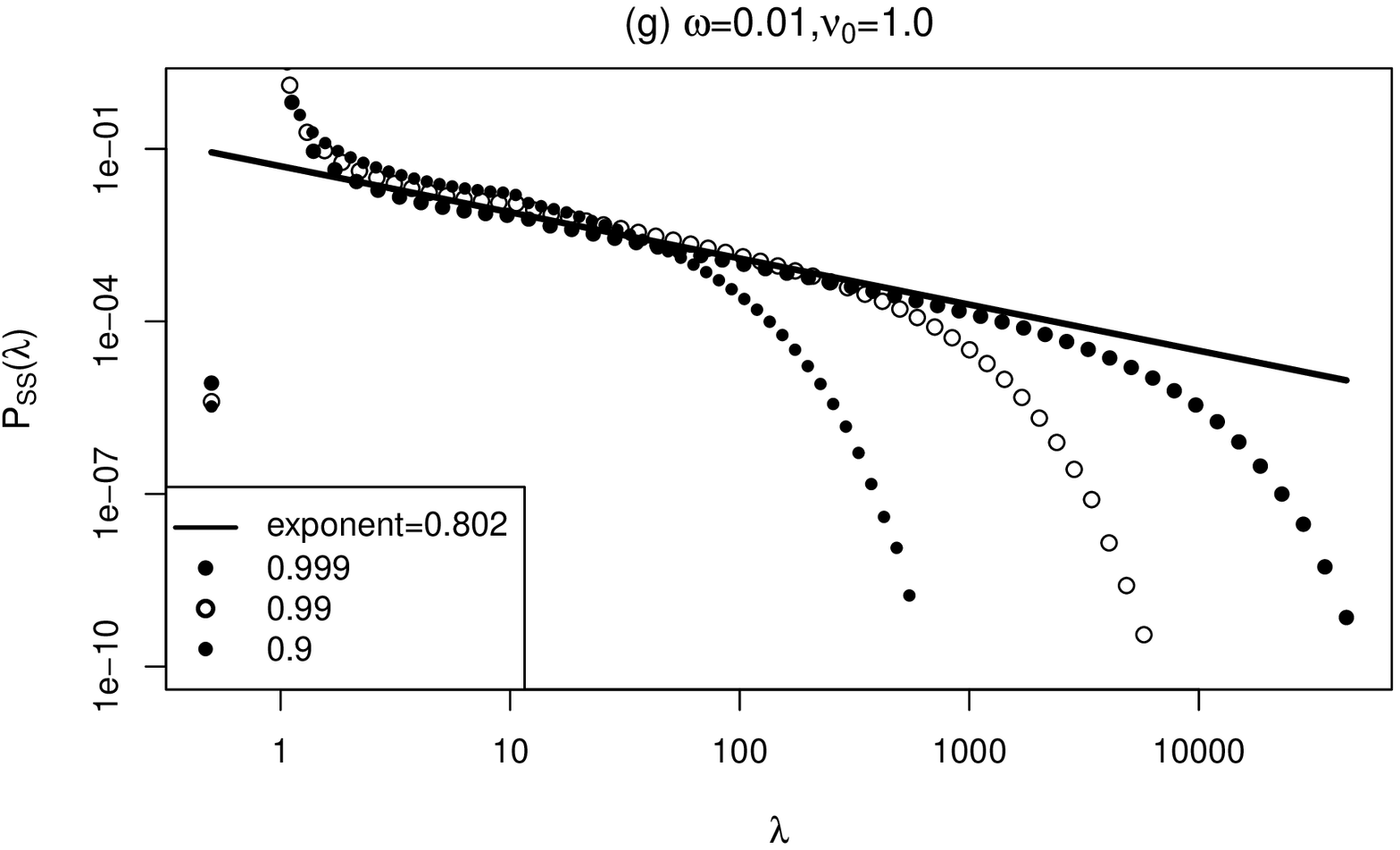}
&   
\includegraphics[width=5cm]{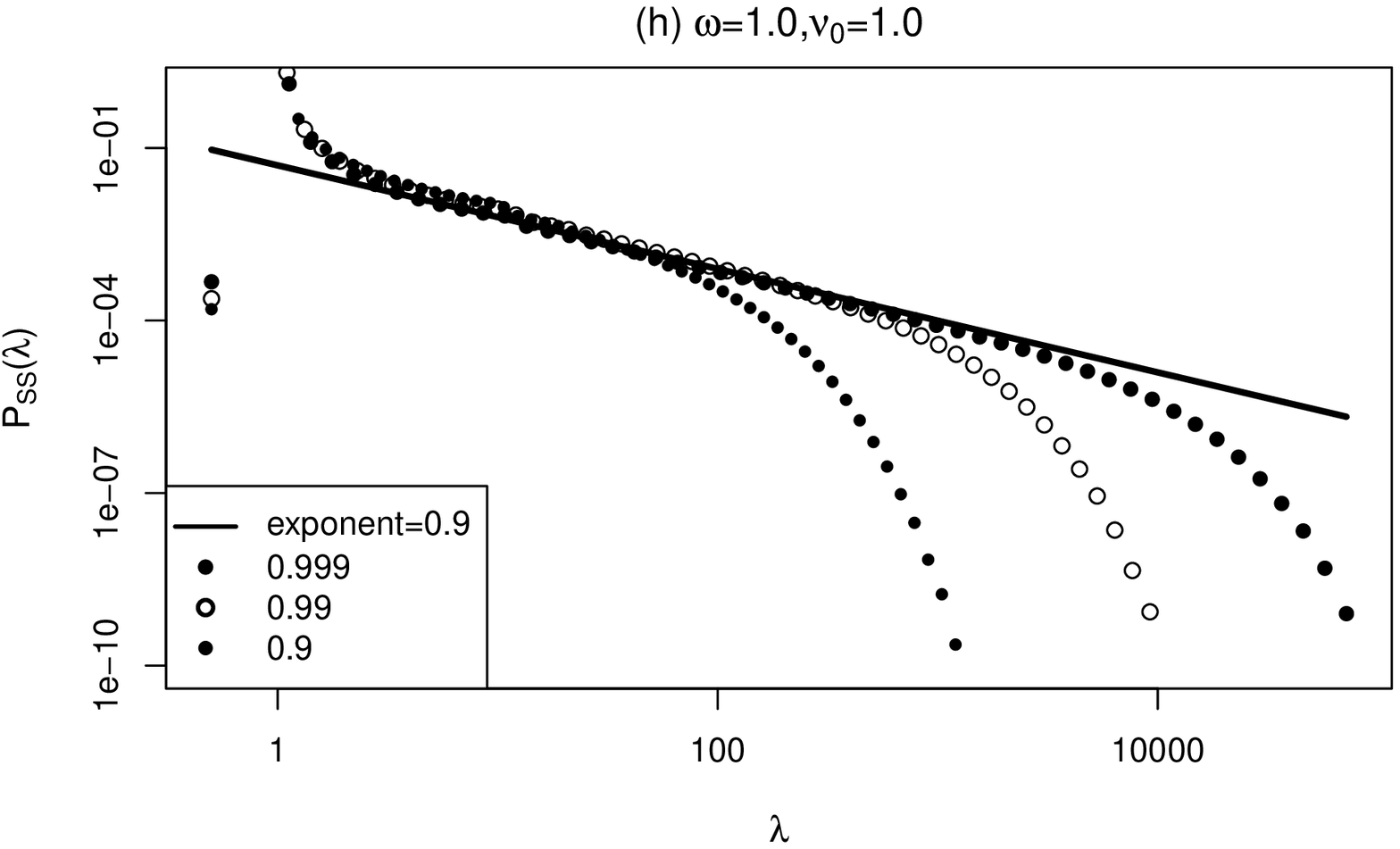}
&   
\includegraphics[width=5cm]{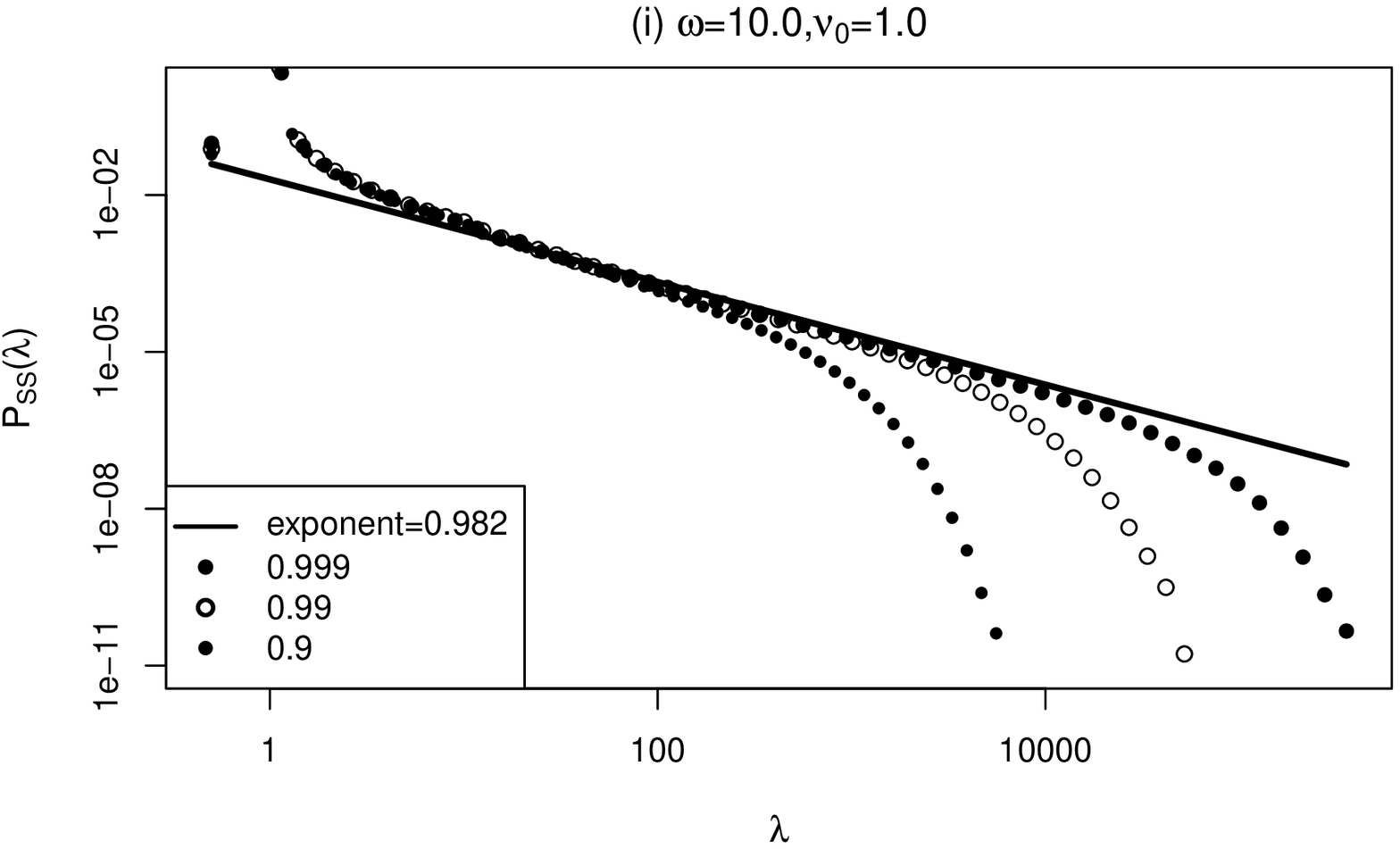}
\end{tabular}
\end{center}
\caption{
  Plot of the numerical results for the steady-state PDF $P_{ss}(\lambda)$ for the general case $K=100$ and $\gamma=11$. 
  We set $n=0.9(\bullet)$, $n=0.99(\circ)$, and $n=0.999(\bullet)$ for $(\tau_1,...,\tau_{100})\sim \mathrm{InvGamma}(11,1),\omega\in\{0.01,1.0,10.0\}$, and $\nu_0\in\{0.01,0.2,1.0\}$.
  (a-c) Background intensity $\nu_0=0.01$ and power-law exponents are $0.99802$, $0.999$, and $0.99982$ for $\omega=0.01$, $1.0$, and $10$, respectively.
  (d-f) Background intensity $\nu_0=0.2$ and power-law exponents are $0.9604$, $0.98$, and $0.9964$ for $\omega=0.01$, $1.0$, and $10$, respectively.
  (g-i) Background intensity $\nu_0=1.0$ and power-law exponents are $0.802$, $0.9$, and $0.982$ for $\omega=0.01$, $1.0$, and $10$, respectively.\color{black}
}
\label{fig:dist_lambda_11}
\end{figure}

We verified the theoretical predictions of Eq.\eqref{eq:power2}
for $\gamma=11$.
\[
P(\lambda)\sim\lambda^{-1+\frac{2\nu_0}{10(\omega+1)}}.
\]
As $\tau\sim \mathrm{InvGamma}(\gamma,1)$ in \eqref{eq:power}, 
we set $K=100$ and $\{\tau_{i}\}_{i=1,\cdots,K}$ as the $(i-1/2)$\% point of the distribution. 
To achieve $n=n_1+...+n_{K}\in \{0.9,0.99,0.999\}$, we set
$n_i=n/K,i=1,\cdots,K$.
Fig.\ref{fig:dist_lambda_11} illustrates the resulting PDFs of 
$\hat{\lambda}(t)$, confirming the validity of the theoretical prediction for $\gamma=11$.

\section{Conclusion}
\label{sec:conclusion}
Here, we investigated the continuous time self-exciting negative binomial process with a memory kernel consisting of a sum of finite $K$ exponential functions.
This process is equivalent to the marked Hawkes process.
Our aim was to extend the previous findings on the intermediate power-law behavior
of the intensities near the critical point for the steady and non-steady state
phase transition to the case of $K>1$.
To achieve this, we developed an efficient sampling method for the marked Hawkes process based on the time-rescaling theorem. We conducted extensive simulations, focusing specifically on the scenario where 
$K=100$ to capture the process 
with a power-law memory kernel. 
Through these simulations, we were able to verify the theoretical predictions and confirm their accuracy.

\appendix

\section{Method of Characteristics}
\label{sec:characteristics}
The method of characteristics is a standard method to solve first-order partial differential equations \citep{gardiner,Kanazawa:2020,Kanazawa:2020-2}.
The equation for the characteristic function 
\[
\phi(s)=\int^{\infty}_{-\infty}\mathrm{e^{isx}}p(x,t|x_0,0)dx
\]
is
\begin{equation*}\label{character-1}
    \partial_t\phi+ks\partial_s\phi=-\frac{1}{2}Ds^2\phi. \tag{A.1}
\end{equation*}
We consider the corresponding Lagrange--Charpit equations
\[
\frac{dt}{dl}=-1, \quad \frac{ds}{dl}=-ks, \quad \frac{d\phi}{dl}=\frac{1}{2}Ds^2\phi
\]
with the parameter $l$ encoding the position along the characteristic curves. 
These equations are equivalent to an invariant form in terms of $l$
\[
 \frac{dt}{1}=\frac{ds}{ks}=-\frac{d\phi}{\frac{1}{2}Ds^2\phi}.
\]
The method of characteristics can be used to solve this equation. Namely, if 
\[
u(s,t,\phi)=a, \qquad \mathrm{and} \qquad  v(s,t,\phi)=b 
\]
are two integrals of the subsidiary equation (with $a$ and $b$ arbitrary constants), then a general solution of \eqref{character-1} is given by 
\[
f(u,v)=0.
\]

\section{Addendum}

In a recent publication by K.Kanazawa and D.Sornette \citep{Kanazawa:2023},
the power-law exponent of the intensity distribution of general
marked Hawkes process was given.
We explain the correspondence for the reader's convenience. 

In \citep{Kanazawa:2023}, the power-law exponent 
is given using the second moment of the mark's PDF E$[m^2]$ as
\[
P_{SS}(\nu)\propto \nu^{-1-a},a=\frac{2\tau\nu_0}{\mbox{E}[m^2]}.
\]
The normalization of $\rho(m)$ with E$[m]=1$ is adopted. 
In our model, $a$ is given as
\[
a=\frac{2\tau\nu_0}{\mbox{E}[m^2]/\mbox{E}[m]}.
\]
We use $\rho(m)$ in Eq.\eqref{PMF} to estimate the power-law exponent.
As E$[m]=\omega/\ln(\omega+1)$, E$[m^2]=\omega(\omega+1)/\ln(\omega+1)$, we obtain 
\[
P_{SS}(\nu)\propto \nu^{-1-a},a=\frac{2\tau \nu_0}{\mbox{E}[m^2]/\mbox{E}[m]}=2\tau\nu_0/(\omega+1).
\]
The result is consistent with ours in \citep{Hisakado:2022}.

\begin{acknowledgements}
This work was supported by JPSJ KAKENHI [Grant No.{} 22K03445].  
We would like to thank Editage (www.editage.com) for English language editing. 
\end{acknowledgements}

%
%

\bibliographystyle{spbasic}      
\bibliography{reference}

\end{document}